\documentclass[11pt]{amsart}
\usepackage{amssymb,amscd,amsthm,amsmath,amsfonts}

\newcommand{\emp}{\emptyset}
\newcommand{\norm}[1]{\mbox{$\left\|#1\right\|$}}
\newcommand{\x}{\times}
\newcommand{\cs}{\mbox{$C^{*}$-algebra}}
\newcommand{\css}{\mbox{$C^{*}$-algebras}}

\newcommand{\ov}[1]{\mbox{$\overline{#1}$}}

\newcommand{\C}{\mathbb{C}}

\newcommand{\Z}{\mathbb{Z}}

\newcommand{\al}{\mbox{$\alpha$}}
\newcommand{\eps}{\mbox{$\epsilon$}}
\newcommand{\bt}{\mbox{$\beta$}}
\newcommand{\ga}{\mbox{$\gamma$}}
\newcommand{\Ga}{\mbox{$\Gamma$}}
\newcommand{\de}{\mbox{$\delta$}}
\newcommand{\De}{\mbox{$\Delta$}}

\newcommand{\la}{\mbox{$\lambda$}}

\newcommand{\Om}{\mbox{$\Omega$}}

\newcommand{\mfB}{\mathfrak{B}}

\newcommand{\mfK}{\mathfrak{K}}
\newcommand{\mfL}{\mathfrak{L}}

\newcommand{\mfH}{\mathfrak{H}}

\newcommand{\bgc}{\begin{center}}

\newcommand{\edc}{\end{center}}
\newcommand{\be}{\begin{enumerate}}
\newcommand{\ee}{\end{enumerate}}
\newcommand{\beqn}{\begin{eqnarray}}
\newcommand{\eeqn}{\end{eqnarray}}
\newcommand{\beqns}{\begin{eqnarray*}}
\newcommand{\eeqns}{\end{eqnarray*}}
\newcommand{\bq}{\begin{quote}}
\newcommand{\eq}{\end{quote}}

\newcommand{\bi}{\begin{itemize}}
\newcommand{\ei}{\end{itemize}}
\newcommand{\bd}{\begin{description}}
\newcommand{\ed}{\end{description}}

\newcommand{\lan}{\mbox{$\langle$}}
\newcommand{\ran}{\mbox{$\rangle$}}

\theoremstyle{plain}
\newtheorem{theorem}{Theorem}

\newtheorem{proposition}{Proposition}
\newtheorem{corollary}{Corollary}

\numberwithin{equation}{section}

\begin{document}
\title{The Fourier algebra for locally compact groupoids}
\author{Alan L. T. Paterson}
\address{Department of Mathematics\\
University of Mississippi\\
University, MS 38677}
\email{mmap@olemiss.edu}
\keywords{Fourier algebra, locally compact groupoids, Hilbert modules,
positive definite functions, completely bounded maps}
\subjclass{43A32}
\date{}
\begin{abstract}
We introduce and investigate using Hilbert modules the properties of the {\em Fourier algebra} $A(G)$ for
a locally compact groupoid $G$.  We establish a duality theorem for such groupoids in terms of multiplicative module maps.  This includes as a special case the classical duality theorem for locally compact groups proved by P. Eymard.
\end{abstract}

\maketitle

\section{Introduction}

For an abelian locally compact group $H$, the Fourier algebra $A(H)$ and the Fourier-Stieljes algebra $B(G)$ are just $L^{1}(\widehat{H})$ and $M(\widehat{H})$ respectively, and are taken by the Fourier transform into certain subalgebras of $C(H)$.
The Fourier algebra $A(G)$ and the Fourier-Stieltjes algebra $B(G)$ for non-abelian locally compact groups $G$ were introduced and studied in the paper of Eymard (\cite{eymard}). These are commutative Banach algebras and are also subalgebras of $C(G)$.  Eymard showed that the character space of $A(G)$ can be identified with $G$ in 
the natural way.  Walter (\cite{wal}) showed that both $A(G)$ and $B(G)$ as Banach algebras determine the group $G$.  Since \cite{eymard,wal}, the
study of $A(G)$ has developed rapidly. 

A corresponding theory of $A(G)$ and $B(G)$
for a locally compact groupoid $G$ has been developed only recently and has
gone in two related directions. The first of these, due to J. Renault, develops
the theory for a measured groupoid $G$.  So a quasi-invariant measure on the unit space is presupposed.  This fits in with the locally compact
group case, the measure on the singleton unit space there being, of course,
just the point mass.  This work has been further developed by 
Jean-Michel Vallin (\cite{vall1,vall}).  Using Hopf-von Neumann bimodule structures, he generalizes Leptin's theorem relating the amenability of the measured groupoid $G$ to the existence of a bounded approximate identity in the Fourier algebra.

The other approach, due to Ramsay and Walter (\cite{ramwal}) starts with a locally compact groupoid without a choice of a quasi-invariant measure.  They
show that there exists a natural candidate for the Fourier-Stieltjes algebra
on $G$, viz. the span $\mathcal{B}(G)$ of the bounded Borel positive definite
functions on $G$. This is then realized as a Banach
algebra of completely bounded bimodule maps using the universal representation
of $G$. Oty (\cite{oty}) investigated the algebra of continuous functions in
$\mathcal{B}(G)$, and, by analogy with the group case, suggested a natural
version of $A(G)$.  

In this paper, we present a continuous version of $A(G)$ parallel to that for the group case. 
As the theory given is a {\em continuous} theory, we concentrate on that part of the 
representation theory of $G$ determined by {\em continuous} $G$-Hilbert bundles over $G^{0}$
rather than that determined by measurable Hilbert bundles for a given quasi-invariant measure.  The canonical example of a
continuous $G$-Hilbert bundle is the Hilbert bundle $L^{2}(G)$ of the $L^{2}(G^{u})$'s
with natural left regular $G$-action, and using $C_{c}(G)$
in the obvious way as a total space of continuous sections for the bundle.

Accordingly, we can formulate an approach to the Fourier-Stieltjes and Fourier
algebras for $G$ as follows. Group representations on a Hilbert space
get replaced in the groupoid case by continuous $G$-Hilbert modules, and the left regular representation in the group case gets replaced in the groupoid case by $L^{2}(G)$.  For a general continuous $G$-Hilbert bundle, we consider the space $\De_{b}$ of continuous bounded sections of the bundle.  This Banach space is a Hilbert $C_{0}(G^{0})$-module.  Of particular importance for this paper is the space $E^{2}$ of continuous bounded sections of $L^{2}(G)$ that vanish at infinity.  Again, $E^{2}$ is a Hilbert $C_{0}(G^{0})$-module.

The Fourier-Stieltjes algebra $B(G)$ is
just the space of ``coefficients'' $(\xi,\eta)$ ($x\to (L_{x}\xi(s(x)),\eta(r(x)))$)
arising from a continuous $G$-Hilbert module with $\xi,\eta\in \De_{b}$. A simple, but important, fact, is that these 
coefficients do not depend on any quasi-invariant measure on $G^{0}$.  
Of course, representation theory comes in when we put a quasi-invariant measure on $G^{0}$ for some $G$-Hilbert bundle. 

In the first approach to norming $B(G)$, we follow the approach of Renault
\cite{ren} to show that $B(G)$ is a Banach algebra, the norm $\norm{\phi}$ of $\phi\in B(G)$ being given by the $\inf$ of 
$\norm{\xi}\norm{\eta}$ over all possible ways of representing $\phi=(\xi,\eta)$.  

We briefly consider another norm on $B(G)$ inspired by the paper \cite{ramwal} of Ramsay and Walters.  In their approach, the norm on $\mathcal{B}(G)$ is defined in terms of the completely bounded multiplier norm on the $\cs$ $M^{*}(G)$.  Here, $M^{*}(G)$ is the completion 
of the image under the universal representation of the convolution algebra of compactly supported, bounded Borel functions on $G$.  Each $\phi\in \mathcal{B}(G)$ acts as a 
multiplier $T_{\phi}$ on $M^{*}(G)$ and this is a completely bounded operator on $M^{*}(G)$.  Define $\norm{\phi}_{cb}=\norm{T_{\phi}}_{cb}$.  They show that $\mathcal{B}(G)$ is a Banach algebra under $\norm{.}_{cb}$.  

In our continuous situation, we follow this approach with $M^{*}(G)$ replaced by $C^{*}(G)$, and we show that $B(G)$ is a normed algebra under the resulting cb-norm $\norm{.}_{cb}$.  Further, 
$\norm{.}_{cb}\leq \norm{.}$ on $B(G)$.  In general, $(B(G),\norm{.}_{cb})$ is not complete (so that $\norm{.}$ is not equivalent to $\norm{.}_{cb}$).  However, by \cite{waldual,paul}, the two norms on $B(G)$ are the same for locally compact groups and for the trivial groupoids 
$G_{n}=\{1,2,\ldots ,n\}\x \{1,2,\ldots ,n\}$.  The two norms always coincide on $P(G)$.

 The Fourier algebra $A(G)$ is defined to be the closure in $B(G)$ of the algebra generated by the coefficients of $E^{2}$ (pointwise operations).  An important subspace $\mathcal{A}(G)$ of $A(G)$ 
is also used in the paper. (The two spaces are the same in the group case.)  We need $\mathcal{A}(G)$ for duality reasons as it relates naturally to $VN(G)$ (below).

In the case of the Hilbert bundle $L^{2}(G)$, $C_{c}(G)$ has natural right and left convolution actions on the Hilbert module $E^{2}$.  The right convolution operators are adjointable and generate the reduced $\cs$ $C^{*}_{red}(G)$ of $G$.  The strong operator closure of the algebra of left convolution operators plays an important role in the theory and is denoted by $VN(G)$.  The operators in $VN(G)$ are rarely adjointable and $VN(G)$ has to be treated in Banach algebra terms.  A useful fact is that $VN(G)$ has a bounded left approximate identity.

With Eymard's theorem in mind, and noting that $A(G)$ is a commutative Banach algebra and that $G$ is contained in the natural way in the set of characters of $A(G)$, we would naturally ask  if the character space of $A(G)$ is equal to $G$.  Eymard's theorem says that this is the case if $G$ is a group.  I do not know if this is true for groupoids.  

We obtain instead another duality result
by appropriately adapting Eymard's group argument to the groupoid case.  This duality result involves sections and $C_{0}(G^{0})$-module maps rather than scalar homomorphisms.  For motivation, consider the case where $G$ is a group bundle $\cup_{u\in G^{0}}G^{u}$.  Each character on $A(G^{u})$ is determined by a point $x^{u}$ of $G^{u}$ and trying to involve the whole of the groupoid $G$ rather than just one $G^{u}$, it is reasonable to think of continuously varying the $x^{u}$'s to get a section $u\to x^{u}$ of $G$.  Such a section determines a multiplicative continuous module map $\ga:A(G)\to C_{0}(G^{0})$.  

For a general groupoid $G$, it is then natural to consider the group $\Ga$ of {\em bisections} of $G$.   A bisection is a subset $a$ of $G$ which determines homeomorphic sections $u\to a^{u}$, $u\to a_{u}$ for the $r-$ and $s-$ maps.  The group $\Ga$ will be our ``dual'' for $G$.  The main theorem of the paper is that for a large class of groupoids $G$, $\Ga$ can be identified with a certain set of multiplicative $C_{0}(G^{0})$-module maps on $A(G)$. Each of the maps determines a pair of left and right module multiplicative maps
linked up by a homeomorphism of $G^{0}$.
We also require a condition on the restriction of these maps to $\mathcal{A}(G)$ which ennables us to use the operators of $VN(G)$ to prove the identification.  We also have to restrict to groupoids for which there are ``many'' bisections and which are locally a product.  (There are many examples of such groupoids.)  

A number of natural open questions are raised throughout the paper, and some of these are listed at the end.  I am grateful to Karla Oty, Arlan Ramsay and Marty Walter for helpful discussions.


\section{Preliminaries}

All locally compact spaces are assumed Hausdorff and second countable, and all Hilbert spaces 
separable.

Let $X$ be a locally compact Hausdorff space. Then $C(X)$ is the space of bounded continuous complex-valued functions on $X$.  The subalgebras of functions that vanish at infinity and that have compact support are respectively denoted by $C_{0}(X)$ and $C_{c}(X)$.  For $f\in C(X)$, 
$S(f)$ denotes the support of $f$.

If $E$ is a Banach space, then $\mfB(E)$ is the Banach algebra of bounded linear operators on $E$.  If $A$ is a $\cs$ and $E,F$ are Hilbert $A$-modules, then $\mfL(E,F), \mfK(E,F)$ 
are the spaces of adjointable and compact maps from $E$ to $F$.  We write $\mfL(E,E)=\mfL(E)$.
Of course $\mfL(E)$ is a $\cs$. A good account of Hilbert $C^{*}$-modules is the book of 
Lance (\cite{lance}).

Throughout the paper, $G$ will stand for a locally compact Hausdorff groupoid. (This class of groupoids is treated in detail in the books \cite{rg,muhlyTCU,pat2} to which the reader is referred for more information.) The unit space of $G$ is denoted by $G^{0}$, and the range and source maps by $r,s$. 
For $u\in G^{0}$,  we define $G^{u}=r^{-1}(\{u\})$ and
$G_{u}=s^{-1}(\{u\})$.  Note that a product $xy$ in $G$ makes sense if and only if 
$s(x)=r(y)$.  We define
$G^{2}=G*G=\{(x,y)\in G\x G: r(y)=s(x)\}$.
The product map $(x,y)\to xy$ is defined on $G^{2}$.  If $A\subset G$, then we write
$A^{u}=A\cap G^{u}, A_{u}=A\cap G_{u}$.  An {\em $r$-section} of $G$ is a subset $A$ of 
$G$ such that for all $u\in G^{0}$, the set $A^{u}$ is a singleton $\{a^{u}\}$, and the map
$u\to a^{u}$ is a homeomorphism from $G^{0}$ onto $A$.  An $s$-section is defined similarly, and a {\em bisection} is a subset $A$ of $G$ that is both an $r$-section and an $s$-section.

As usual, a left Haar system $u\to \la^{u}$ is presupposed on the fibers $G^{u}$, and 
$\la_{u}=(\la^{u})^{-1}$ (a measure on $G_{u}$).  
For the sake of brevity, we write $D=C_{0}(G^{0})$. Useful norms on $C_{c}(G)$ are the $I$-norm's, 
$\norm{.}_{I,r}, \norm{.}_{I,s}$ and $\norm{.}_{I}$.  Here, for $f\in C_{c}(G)$,
\[ \norm{f}_{I,r} = \sup_{u\in G^{0}}\int_{G^{u}}\mid f(t)\mid\,d\la^{u}(t), \]
\[ \norm{f}_{I,d} = \sup_{u\in G^{0}}\int_{G_{u}}\mid f(t)\mid\,d\la_{u}(t)\]
and
\[ \norm{f}_{I} = \max\{\norm{f}_{I,r},\norm{f}_{I,d}\}.  \]
Note that $(C_{c}(G), \norm{.}_{I,r})$ is a normed algebra, whereas $(C_{c}(G),\norm{.}_{I})$ is a normed $^{*}$-algebra (isometric involution).  The involution on $C_{c}(G)$ is the map
$f\to f^{*}$ where $f^{*}(x)=\ov{f(x^{-1})}$.  We also define $f^{\vee}\in C_{c}(G)$ by: $f^{\vee}(x)=f(x^{-1})$.

We next recall some details concerning the disintegration of representations of $C_{c}(G)$.  The theorem is due to J. Renault (\cite{renrep}).  A detailed account of the theorem is given in the book of Paul Muhly (\cite{muhlyTCU}).  Let $\Phi$ be a representation of $C_{c}(G)$ on a Hilbert space $H$ which is continuous in the inductive limit topology.  Then $\Phi$ disintegrates as follows.  

There is a probability measure $\mu$ on $G^{0}$ which is {\em quasi-invariant} (in the sense defined below).  Associated with $\mu$ is a positive regular Borel measure $\nu$ on $G$ defined by: $\nu=\int\la^{u}\,d\mu(u)$.
(When we want to make the dependence of $\nu$ on $\mu$ explicit, we will use (as in \cite{ramwal}) $\la^{\mu}$ in place of $\nu$.)
The measure $\nu^{-1}$ is the image under $\nu$ by inversion: precisely, $\nu^{-1}(E)=\nu(E^{-1})$.  There is also a measure $\nu^{2}$ on $G^{2}$ given by: $\nu^{2}=\iint \la_{u}\x \la^{u}\,d\mu(u)$. The quasi-invariance of $\mu$ just means that $\nu$ is equivalent to $\nu^{-1}$.  The {\em modular function} $D$ is defined as the Radon-Nikodym derivative $d\nu/d\nu^{-1}$. The function $D$ can be taken to be Borel (\cite{hahn,ramsayJFA,muhlyTCU}), and satisfies the 
properties: $D(x^{-1})=D(x)^{-1}$ $\nu$-almost everywhere, and 
$D(xy)=D(x)D(y)$ $\nu^{2}$-almost everywhere.  Let $\nu_{0}$ be the measure on $G$ given by: $d\nu_{0}=D^{-1/2}d\nu$.

Next, there exists a $\mu$-measurable Hilbert bundle $\mfK$ over $G^{0}$ and a {\em $G$-representation} $L$ on $\mfK$.  So each $L(x)$ ($x\in G$) is a linear isometry from $\mfK_{s(x)}$ onto $\mfK_{r(x)}$, which is multiplicative $\nu^{2}$-almost everywhere and inverse preserving 
$\nu$-almost everywhere.  Further, for every pair of measurable sections $\xi,\eta$ of $\mfK$, it is required that the function $x\to \lan L(x)\xi(s(x)),\eta(r(x))\ran$ is $\mu$-measurable.  The representation $\Phi$ of $C_{c}(G)$ is then given by:
\begin{equation}  
\lan \Phi(F)\xi,\eta\ran = \int F(x)\lan L(x)\xi(s(x)),\eta(r(x))\ran\,d\nu_{0}(x).   \label{eq:disinteg}
\end{equation}
We will refer to the triple $(\mu,\mfK,L)$ as a {\em representation} of $G$.

Of particular importance is the regular representation $\pi_{u}$ for $u\in G^{0}$.  Here, for 
$F\in C_{c}(G)$, $\pi_{u}:C_{c}(G)\to \mfB(L^{2}(G_{u}))$ is the representation given by:
\begin{equation}     \label{eq:reg}
\pi_{u}(F)(\xi)(x)=\int F(t)\xi(t^{-1}x)\,d\la^{r(x)}(t)=F*\xi(x).
\end{equation}
(The quasi-invariant measure and measurable bundle associated with $\pi_{u}$ is calculated in 
\cite[Example 3.26]{muhlyTCU}.)  We define a $C^{*}$-norm $\norm{.}_{red}$ on $C_{c}(G)$ by setting $\norm{F}_{red}=\sup_{u\in G^{0}}\norm{\pi_{u}(F)}$.  The reduced $\cs$ 
$C_{red}^{*}(G)$ is defined to be the completion of $(C_{c}(G),\norm{.}_{red})$.

We conclude this section with the following
simple proposition.  This is surely well-known - a smooth version of it for pseudodifferential operators is given in \cite{peetre} - but for the convenience of the reader, we give a proof.

\begin{proposition}        \label{prop:supp}
Let $X$ be a locally compact Hausdorff space and $R:C_{0}(X)\to C_{0}(X)$ be a bounded linear map such that for all $f\in C_{c}(X)$, we have $S(Rf)\subset S(f)$.  Then there exists a bounded continuous function $k:X\to \C$ such that $Rf=kf$ for all $f\in C_{0}(X)$.
\end{proposition}
\begin{proof}
Let $f\in C_{c}(X)$, $x_{0}\in X$ and $\phi\in C_{c}(X)$ be such that $\phi=1$ on a neighborhood of $x_{0}$. Let $g=f-f(x_{0})\phi$.  Then $g(x_{0})=0$.  There exists a sequence $\{g_{n}\}$ in 
$C_{c}(X)$ such that $g_{n}\to g$ and $g_{n}=0$ on a neighborhood $U_{n}$ of $x_{0}$.  Since
$S(Rg_{n})\subset S(g_{n})$, we have $S(Rg_{n})\subset X\setminus U_{n}$.  So $Rg_{n}(x_{0})=0$.  Since $Rg_{n}\to Rg$, we have $Rg(x_{0})=0$.  So 
$Rf(x_{0})=R(g)(x_{0})+f(x_{0})R\phi(x_{0})=((R\phi)f)(x_{0})$.  Now if $\psi\in C_{c}(X)$ is such that $\psi=1$ on neighborhood of $x_{0}$, then $\phi-\psi$ vanishes on a neighborhood of $x_{0}$ and $R(\phi-\psi)(x_{0})=0$.  So we obtain a well-defined function $k$ on $X$ by setting $k(x_{0})=R\phi(x_{0})$.  Clearly, $k$ is continuous since the same $\phi$ applies on a neighborhood of $x_{0}$, and $\norm{k}\leq \norm{R}$ (since we can arrange that $\norm{\phi}=1$).  Lastly, $Rf=kf$ for $f\in C_{c}(X)$ and by continuity, this is valid for 
$f\in C_{0}(X)$.
\end{proof}     


\section{A Fourier-Stieltjes algebra $B(G)$}

In this section, we discuss the basic facts about the (continuous) Fourier-Stieltjes algebra $B(G)$ for a locally compact groupoid.  The approach and proofs are inspired by those of the papers \cite{ren} by Jean Renault for measured groupoids and 
\cite{ramwal} by Arlan Ramsay and Marty Walters for the Borel case.
In order to fit in with the customary practice of the theory of Hilbert $C^{*}$-modules, we will assume that the   
Hilbert modules (and spaces) are conjugate linear in the first variable.

Let $X$ be a locally compact Hausdorff space and 
$\mfH=\{\mfH^{u}\}$ $(u\in X)$ be a {\em continuous} Hilbert bundle over
$X$. In the notation of Dixmier (\cite[Ch.10]{D2}), $\mfH$ is just a
continuous field of Hilbert spaces over $X$.   The norm on each $\mfH^{u}$ is denoted by 
$\norm{.}^{u}$.  Let $\De_{b}$ be the space of continuous bounded sections of $\mfH$.  Then 
$\De_{b}$ is a Banach space with norm given by $\norm{\xi}=\sup_{u\in X}\norm{\xi(u)}^{u}$.  The space $\De_{0}$ of sections $\xi\in \De_{b}$ of $\mfH$ that vanish at infinity, i.e. such that $\norm{\xi}^{u}\to 0$ as $u\to \infty$ in $X$, is a closed subspace of $\De_{b}$.  The space $\De_{c}$ of elements of $\De_{0}$ with compact support on $X$ is a dense subspace of $\De_{0}$.  
We will write $\De_{b}(\mfH), \De_{0}(\mfH)$ and $\De_{c}(\mfH)$ instead of $\De_{b}, \De_{0}$ and $\De_{c}$ when we wish to make explicit $\mfH$.
By \cite[Proposition 10.1.9]{D2}, all of the spaces $\De_{b}, \De_{0}, \De_{c}$ are $C_{0}(X)$ modules under the action: $(a,\xi)\to a\xi$, where $a\xi(u)=a(u)\xi(u)$.  
Note that if $u\in G^{0}$ then the linear space $\{\xi(u): \xi\in \De_{c}\}$ is dense in $\mfH^{u}$.

Now let $G$ be a locally compact groupoid.  
A Hilbert bundle  $\mfH$ over $G^{0}$ is called a {\em $G$-Hilbert bundle} if for each $x\in G$, there is given a linear isometry $L_{x}$ from $H_{s(x)}$ onto $H_{r(x)}$ such for each $\xi,\eta\in \De_{b}$, the map
$x\to (L_{x}\xi(s(x)),\eta(r(x)))$ is continuous, and the map $x\to L_{x}$ is a groupoid homomorphism from $G$ into the isomorphism groupoid of the fibered set $\cup_{u\in G^{0}}\mfH^{u}$ (\cite[Ch. 1]{muhlyTCU}).   As in \cite{ren}, we will denote the function
$x\to (L_{x}\xi(s(x)),\eta(r(x)))$ on $G$ by $(\xi,\eta)$, and will call $(\xi,\eta)$
a {\em coefficient} of the Hilbert bundle $\mfH$.

\begin{proposition}     \label{prop:inftyxieta}
$\norm{(\xi,\eta)}_{\infty}\leq \norm{\xi}\norm{\eta}$.
\end{proposition}
\begin{proof}
$\mid (\xi,\eta)(x)\mid=\mid (L_{x}\xi(s(x)),\eta(r(x)))\mid
\leq \norm{\xi(s(x))}^{s(x)}\norm{\eta(r(x))}^{r(x)}\leq \norm{\xi}\norm{\eta}$.
\end{proof}

A function $\phi\in C(G)$ is said (\cite{ramwal}) to be {\em positive definite} if for all $u\in G^{0}$ and all $f\in C_{c}(G)$ we have 
\begin{equation}
\iint \phi(y^{-1}x)f(y)\ov{f(x)}\,d\la^{u}(x)\,d\la^{u}(y)\geq 0.  \label{eq:posd}
\end{equation}
The set of all positive definite functions in $C(G)$ is denoted (\cite{oty}) by $P(G)$.
It is easy to check that if $\xi\in \De_{b}(\mfH)$, then $(\xi,\xi)\in P(G)$.   The converse is also 
true.

\begin{theorem}      \label{th:contpos}
Let $\phi\in C(G)$.  Then $\phi$ is positive definite if and only if $\phi$ is a coefficient
of the form $(\xi,\xi)$ for some $G$-Hilbert bundle.
\end{theorem}
\begin{proof}
The Borel version of this theorem is proved in \cite[Theorem 3.5]{ramwal} and the proof in the
continuous case is easier.  The relevant $G$-Hilbert bundle $\mathfrak{H}$ is that determined by the semi-inner product spaces $L^{2}(G^{u})$ where 
\[   (f\mid g)_{u}=\iint \phi(y^{-1}x)g(y)\ov{f(x)}\,d\la^{u}(x)\,d\la^{u}(y).  \]
This gives a continuous Hilbert bundle (since the $\la^{u}$'s vary ``continuously'').
The continuity of the section $\xi$ is given in \cite[Theorem 3.5]{ramwal}.
\end{proof}

\begin{proposition}           \label{prop:posd2}
A function $\phi\in C(G)$ is positive definite if and only if for any $n$, any $u\in G^{0}$,
any $x_{1},\ldots ,x_{n}\in G^{u}$ and any $\al_{1},\ldots ,\al_{n}\in \C$, we have
\begin{equation}      \label{eq:posd2}
\sum_{i,j}\al_{i}\ov{\al_{j}}\phi(x_{i}^{-1}x_{j})\geq 0.
\end{equation}
\end{proposition}
\begin{proof}
If $\phi\in P(G)$, then direct checking (using Theorem~\ref{th:contpos}) shows that $\phi$ satisfies (\ref{eq:posd2}).  Conversely, approximating $f\,\la^{u}$ weak$^{*}$ by measures of the form 
$(\sum_{i}\al_{i}\de_{x_{i}})$ ($r(x_{i})=u$) gives (\ref{eq:posd}).
\end{proof}

Let $\mfH, \mfK$ be $G$-Hilbert bundles with groupoid actions $x\to L_{x}$, $x\to L'_{x}$ respectively.  Clearly, the direct sums and tensor products of $G$-Hilbert bundles
are themselves are $G$-Hilbert bundles in the natural way.  For example, the fiber Hilbert spaces of 
$\mfH\otimes \mfK$ are $H^{u}\otimes K^{u}$, and the continuous bundle structure is determined by the linear span of sections of the form $\xi\otimes \eta$ where $\xi\in \De_{b}(\mfH), 
\eta\in \De_{b}(\mfK)$.  The groupoid action is given by $x\to L_{x}\otimes L'_{x}$.  

Every fixed Hilbert space $H$ gives rise to a {\em trivial} $G$-Hilbert bundle as follows.  Define $\mfH^{u}=H$ for all $u\in G^{0}$ and take $\De_{b}=C(G^{0},H)$.  The $G$-action on
$\mfH=G^{0}\x H$ is given by: $L_{x}(s(x),z)=(r(x),z)$.  

The set of functions $(\xi,\eta)$ with $\xi,\eta\in \De_{b}(\mfH)$ for some $G$-Hilbert bundle 
$\mfH$ is denoted by $B(G)$ and is called the {\em Fourier-Stieltjes algebra} of $G$.  (Our 
$B(G)$ is the same as the $B_{1}(G)$ of \cite{oty}.)  From the existence of direct sums and tensor products of $G$-Hilbert bundles, it follows that $B(G)$ is an algebra under pointwise operations.  By Theorem~\ref{th:contpos}, $P(G)\subset B(G)$.
Using the polarization identity, every element of $B(G)$ is a linear combination of elements of $P(G)$.  

The norm $\norm{.}$ of $\phi\in B(G)$ is defined by:
\[     \norm{\phi}=\inf\norm{\xi}\norm{\eta}               \]
the $\inf$ being taken over all representations $\phi=(\xi,\eta)$.  

\begin{proposition}        \label{prop:posu}
If $\phi\in B(G)$, then $\norm{\phi}_{\infty}\leq \norm{\phi}$.
If $\phi\in P(G)$, then $\norm{\phi}=\norm{\phi_{\mid G^{0}}}_{\infty}$.  
\end{proposition}
\begin{proof}
The first part follows from Proposition~\ref{prop:inftyxieta}.  For the second, 
write $\phi=(\xi,\xi)$.  For $u\in G^{0}$, $0\leq \phi(u) = \norm{\xi(u)}^{2}$, and so
\[\norm{\xi}^{2}=\norm{\phi_{\mid G^{0}}}_{\infty}\leq \norm{\phi}\leq \norm{\xi}^{2}.  \]
\end{proof}

The proof that $B(G)$ is a Banach algebra relies, as in the corresponding results of
\cite{ren,ramwal}, on a result of \cite{ren} ennabling one to estimate $B(G)$ norms by relating $\phi\in B(G)$ to an element $F\in P(G\x I_{2})$ where $I_{2}$ is the trival groupoid $\{1,2\}\x\{1,2\}$.  This can be regarded as a groupoid version of Paulsen's ``off-diagonalization technique'' (\cite[Th. 7.3]{paul}, \cite[Th. 5.3.2]{effrosruan}).
There seems to be a slight gap in the proof of this result and so we 
give a complete proof of the result for our situation.

The elements of $G\x I_{2}$ are of the form $(x,i,j)$ with $x\in G, i,j\in \{1,2\}$.  We
identify $(G\x I_{2})^{0}$ with $G^{0}\x \{1,2\}$.
A function $F:G\x I_{2}\to \C$ will be identified in the natural way
with the $2\x 2$ matrix-valued function  
\begin{equation*}
\begin{bmatrix}
F(x11) & F(x12)  \\
F(x21) & F(x22)
\end{bmatrix}
\end{equation*}
on $G$.

\begin{proposition}           \label{prop:bgpos}
$\phi\in B(G)$ if and only if there exist $F\in P(G\x I_{2})$ of the form
\begin{equation} \label{eq:bmatr}
\begin{bmatrix}
\rho(x) & \phi(x)  \\
\phi^{*}(x) & \tau(x)
\end{bmatrix}
\end{equation}
with $\rho, \tau\in P(G)$.
\end{proposition}
\begin{proof}
It is proved in \cite[Proposition 1.3]{ren} that given $\phi\in B(G)$, there exists an $F$ of the required form.  Conversely, as in \cite[Proposition 1.3]{ren}, given such an $F$, there exists a Hilbert $(G\x I_{2})$-bundle $\mfH$, with $G\x I_{2}$-action $xij\to L(xij)$,  
and a section $\zeta=(\zeta_{1},\zeta_{2})$ of $\mfH$ such that 
$F=(\zeta,\zeta)$.  Let $\mfH_{i}$ be the restriction of $\mfH$ to $G^{0}\x \{i\}$.  The 
$\mfH_{i}$ are Hilbert $G$-bundles in the obvious way.  The direct sum $\mfH'$ of the 
$\mfH_{i}$'s is then a Hilbert $G$-bundle.  For each $x\in G$, define 
$L'(x):\mfH'_{s(x)}\to \mfH'_{r(x)}$ by:
\[L'(x)=1/2 
\begin{bmatrix}
L(x11) & L(x12)  \\
L(x21) & L(x22)
\end{bmatrix}
\]

It is true that $L'$ is continuous and $L'(xy)=L'(x)L'(y), L'(x)^{*}=L'(x^{-1})$.  However, $L'$ is not usually invertible.  (For example, if $G$ has one point $\{e\}$, $\mfH=\C^{2}$ with the obvious $I_{2}$-action, then $L'(e)$ is the (singular) scalar two-by-two matrix whose entries are all $1$.)  To deal with this, we ``cut down'' $\mfH'$ as follows.  For each $u\in G^{0}$, let $P_{u}$ be the projection $L'(u)$ on $\mfH'_{u}$ and 
$\mfK_{u}=P_{u}(\mfH'_{u})$.  Then $\mfK=\{\mfK_{u}\}$ is a Hilbert bundle in a natural way.  Indeed, let $Y$ be the sections of $\mfH'$ of the form $\xi=\xi_{1}\oplus \xi_{2}$ where $\xi_{i}$ is a continuous section of $\mfH_{i}$.  Then let
$X$ be the vector space of sections of $\mfK$ of the form $P\xi: u\to P_{u}\xi(u)$ where 
$\xi\in Y$.  Then $X$ satisfies (ii) and (iii) of 
\cite[10.1.2]{D2} using the fact that the functions $u\to (P_{u}\xi(u),\eta(u))$ are continuous.  So (\cite[10.2.3]{D2}) $\mfK$ is a Hilbert bundle.  For $x\in G$, let 
$M(x)=L'(x)P_{s(x)}$.  It is left to the reader to check that 
$\mfK$ is a $G$-Hilbert bundle with $G$-action given by $M$.  Then $\phi =2(\eta,\xi)$ where
$\xi=P(\zeta_{1},0), \eta=P(0,\zeta_{2})$, and $\phi\in B(G)$.  
\end{proof}

\begin{theorem}       \label{th:bg}
The set $B(G)$ is a unital commutative Banach algebra under pointwise operations on $G$.
\end{theorem}
\begin{proof}
The proof is effectively the same as the corresponding result (\cite[Proposition 1.4]{ren}) in the measured groupoid case, with $L^{\infty}(G)$ being replaced by $C(G)$ and measurable Hilbert bundles replaced by continuous Hilbert bundles.  That $B(G)$ is unital follows using the trivial
$G$-Hilbert bundle $G^{0}\x \C$ (or by using Proposition~\ref{prop:posd2}).
\end{proof}       

\begin{proposition}       \label{prop:eleag}
\bi
\item[(i)] If $\phi\in B(G)$, then both $\ov{\phi}, \phi^{*}\in B(G)$, and 
$\norm{\phi}=\norm{\ov{\phi}}=\norm{\phi^{*}}$.  
\item[(ii)] $C(G)$ is a two-sided $D$-module, where for $b\in D$ and $F\in C(G)$, the action is given by:
\begin{equation}  \label{eq:bphib}
  (Fb)(x)=F(x)b(r(x)),  (bF)(x)=b(s(x))F(x).  
\end{equation}  
Further, if $f,g\in C_{c}(G)$, then  
\begin{equation}       \label{eq:baction}
 (f*g)b=fb*g, b(f*g)=f*bg.   
\end{equation}
and if $\phi=(\xi,\eta)\in B(G)$, then $\phi b=(\xi,b\eta), b\phi=(\ov{b}\xi,\eta)$.  
Lastly, $B(G)$ is a $D$-submodule of $C(G)$,
and for $\phi\in B(G), b\in D$, we have  
\begin{equation} \label{eq:normphib}
\norm{\phi b}\leq \norm{\phi}\norm{b}, \norm{b\phi}\leq \norm{b}\norm{\phi}.
\end{equation}
\ei
\end{proposition}
\begin{proof}
(i) From the definition of $P(G)$, it follows that 
$\ov{\phi}, \phi^{*}\in P(G)$ whenever $\phi\in P(G)$. Let $\phi\in B(G)$.
Since $B(G)$ is the span of $P(G)$, it follows that $\ov{\phi}, \phi^{*}\in B(G)$.  Let 
$\phi=(\xi,\eta)$ for some $G$-Hilbert bundle $(\mfH,L)$. 

As in the group case, there is a conjugate $G$-Hilbert bundle $(\{\ov{H_{u}}\},\ov{L})$ for any given $G$-Hilbert bundle $(\{H_{u}\},L)$.  Here, if $H$ is  Hilbert space, then $\ov{H}$ is the Hilbert space that coincides with $H$ except that the inner product $( ,)'$ and scalar multiplication $(a,\xi)\to
a.\xi$ for $\ov{H}$ are given by: $(\xi,\eta)'=(\eta,\xi)$ and $a.\xi=\ov{a}\xi$.  We take 
$\ov{L}=L$.  One readily checks that $\ov{\phi}(x)=(\ov{L}_{x}\xi,\eta)'$, so that 
$\norm{\ov{\phi}}\leq \norm{\xi}\norm{\eta}$.  It follows that $\norm{\phi}=\norm{\ov{\phi}}$.
Since $\phi^{*}=(\eta,\xi)$, we obtain $\norm{\phi}=\norm{\phi^{*}}$.

For (ii), see \cite[p.59]{rg}. (Left and right actions of $D$ are interchanged from those in \cite{rg} for duality reasons.) For example, for (\ref{eq:baction}), 
\[ (f*g)b(x)=b(r(x))\int f(t)g(t^{-1}x)\,d\la^{r(x)}(t)=fb*g(x).  \]
\end{proof}

We now turn to the other way of norming $B(G)$. The norm is a
completely bounded norm (\cite{effrosruan}) and the result is a variation of \cite[Theorem 6.1]{ramwal} which says that $\mathcal{B}(G)$
is a Banach algebra. Let $\pi:C_{c}(G)\to C^{*}(G)$ be the canonical isomorphism, and for
$\phi\in B(G)$, define a map $T_{\phi}:\pi(C_{c}(G))\to \pi(C_{c}(G))$ by: 
$T_{\phi}\pi(f)=\pi(\phi f)$.

\begin{theorem}                \label{th:agnorm}
If $\phi=(\xi,\eta)\in B(G)$ is a coefficient of a continuous Hilbert bundle $(\{\mfH_{u}\},L)$, then $T_{\phi}$ is completely bounded on $C^{*}(G)$, and  
\begin{equation}  \label{eq:phifg}
\norm{T_{\phi}}_{cb}\leq \norm{\xi}\norm{\eta}.
\end{equation}
\end{theorem}
\begin{proof} Let $(\mu',\{\mfH_{u}'\},L')$ be a representation of $G$.
Let $\pi'$ be the representation of $C_{c}(G)$ obtained by integrating this representation.  Then $(\mu',\{\mfH_{u}'\otimes \mfH_{u}\},L'\otimes L)$ is a $G$-representation.  Let
$\Tilde{\pi}$ be its integrated form.
Let $n\geq 1$, $A=[f_{ij}]\in M_{n}(C_{c}(G)), A\phi=[f_{ij}\phi]$ and 
$\xi_{i}',\eta_{i}'$ $(1/\leq i\leq n)$ be square integrable sections of $\{\mfH'_{u}\}$.  Let $\xi'=[\xi_{i}]', \eta'=[\eta_{i}]'$.  To obtain the complete boundedness of $T_{\phi}$ and
(\ref{eq:phifg}), it is sufficient to show that 
\begin{equation}              \label{eq:toshow}
\mid (\pi'(A\phi)(\xi'),\eta')\mid\leq \norm{\xi}\norm{\eta}\norm{\xi'}\norm{\eta'}\norm{\Tilde{\pi}(A)}_{n}.   
\end{equation}  
Indeed
\begin{align*}
(\pi'(A\phi)(\xi'),\eta') & =
\sum \int f_{ij}(t)(L_{t}'\xi_{j}(s(t)),\eta_{i}(r(t)))\phi(t)\,d\nu_{0}(t)\\
& =\sum\int f_{ij}((L_{t}'\otimes L_{t})((\xi_{j}'\otimes \xi)(s(t)),
(\eta_{i}'\otimes \eta)(r(t))\,d\nu_{0}(t)\\
&=(\Tilde{\pi}(A)(\xi'\otimes \xi),\eta'\otimes \eta).  
\end{align*}
So  
$\mid (\pi'(A\phi)\xi',\eta')\mid\leq \norm{\pi(A)}_{n}\norm{\xi}\norm{\eta}\norm{\xi'}\norm{\eta'}$ and 
(\ref{eq:toshow}) follows.
\end{proof}
\begin{corollary}   \label{cor:norms}
$\norm{.}_{cb}\leq \norm{.}$ on $B(G)$.
\end{corollary}

In the Borel case, Ramsay and Walter (\cite[Theorem 6.1]{ramwal}) show that $\mathcal{B}(G)$
is a Banach algebra.  However, $(B(G),\norm{.}_{cb})$ is not always a Banach algebra.  This follows from the elegant counterexample at the end of \cite{ramwal} for which $G$ is a bundle of groups $X\x \Z$ with $X=\{re^{\imath\theta}: 0\leq r\leq 1,\theta\in \{0,1,1/2,1/3,\ldots \}\}$.
In some cases, however, $\norm{.}=\norm{.}_{cb}$ on $B(G)$.  By \cite{waldual,paul}, this is the case for locally compact groups and also for the case of the trivial groupoids $G_{n}=\{1,2,\ldots ,n\}\x \{1,2,\ldots ,n\}$.  (For the $G_{n}$ case, see \S 5, Example below.)
The measured groupoid version of the result has been proved in complete generality by Renault
(\cite[Theorem 22]{ren}) using the module Haagerup tensor product.  Note that from 
Corollary~\ref{cor:norms} and Banach's isomorphism theorem, if $(B(G),\norm{.}_{cb})$ is a Banach algebra, then the norms $\norm{.}, \norm{.}_{cb}$ are equivalent on $B(G)$.

\section{The left regular Hilbert bundle}

Let $L^{2}(G)=\{L^{2}(G^{u})\}$.  In the natural way, this is a $G$-Hilbert bundle, which we will call the {\em left regular} Hilbert bundle of $G$.  (This Hilbert bundle has been used by Khoshkam and Skandalis (\cite{khoshskand}) even in the non-Hausdorff context.)
In more detail, we regard the functions $f\in C_{c}(G)$ in the obvious way as sections $u\to f^{u}=f_{\mid G^{u}}$ which determine the continuous sections of $L^{2}(G)$.   In fact, the space of such sections satisfies axioms (i), (ii) and (iii) of \cite[Definition 10.1.2]{D2}, and so (\cite[Proposition 10.2.3]{D2}) determines a continuous field, the continuous sections in general just being those sections that are locally close to $C_{c}(G)$.  The $G$-action on 
$\De_{b}(L^{2}(G))$ is given by: $(L_{x}\xi)(t)=\xi(x^{-1}t)$.  

In particular, a section $u\to F^{u}$ of $L^{2}(G)$ is continuous iff the map $u\to \norm{F^{u}}$ is continuous and for all $g\in C_{c}(G)$, the function $u\to \lan F^{u},g^{u}\ran$ is continuous for all $g\in C_{c}(G)$.  Let $E^{2}(G)$, or simply $E^{2}$, be the set of continuous sections of $L^{2}(G)$ that vanish at infinity. Of course, $C_{c}(G)\subset E^{2}$, and $E^{2}$ is a Banach space under the section norm: $\norm{F}=\sup_{u\in G^{0}}\norm{F^{u}}$.  Also in the canonical way, $E^{2}$ is a $D$-module: for $\xi\in E^{2}$ and $b\in D$, we set
$b\xi(t)=(\xi b)(t)=\xi(t)b(r(t))$.

\begin{proposition}         \label{prop:ccgdense}
$C_{c}(G)$ is dense in $E^{2}$.
\end{proposition}
\begin{proof}
Let $F\in E^{2}, \eps>0$.  We show that there exists $f\in C_{c}(G)$ such that 
$\norm{F-f}<\eps$.  Since $F$ vanishes at infinity, we can suppose (by multiplying $F$ by a suitable $b\in D$) that $F$ has compact support $C\subset G^{0}$.  Let 
$u\in C$ and $f[u]\in C_{c}(G^{u})$ be such that 
$\norm{F^{u}-f[u]}<\eps$.  Extend $f[u]$ to a function $f[u]'\in C_{c}(G)$.  Then there exists a neighborhood $W(u)$ of $u$ in $G^{0}$ such that $\norm{F^{v}-(f[u]')^{v}}<\eps$ for all $v\in W(u)$.  Cover $C$ by a finite number of the $W(u)$'s, say $W(u_{1}),\ldots ,W(u_{n})$ and let $\{b_{i}\}$ be 
such that $b_{i}\in C_{c}(W(u_{i}))$ and $b_{i}\geq 0$, $\sum_{i=1}^{n}b_{i}=1$ on 
$C$.  Then $\sum_{i=1}^{n}b_{i}F=F$, and 
$\norm{F-\sum_{i=1}^{n}b_{i}f[u_{i}]'}<\eps$.              
\end{proof}

If $G$ is a locally compact group, then the reduced $\cs$ $C_{red}^{*}(G)$ and its enveloping
von Neumann algebra $VN(G)$ are defined on the Hilbert space $L^{2}(G)$.  We need versions of these for the groupoid case.  In the groupoid case, $L^{2}(G)$ is replaced by $E^{2}$.  While $E^{2}$ is not a Hilbert space, it is a Hilbert $D$-module.  The right $D$-action has been given above, while the $D$-valued inner product $\lan ,\ran$ on $E^{2}$ is given by:
\[     \lan \xi,\eta\ran(u)=(\xi^{u},\eta^{u}).   \]

In this Hilbert module context, $C^{*}_{red}(G)$ will be the $C^{*}$-subalgebra of $\mfL(E^{2})$
generated by the right regular antirepresentation of $C_{c}(G)$.  We will take $VN(G)$ to be the commutant of $C^{*}_{red}(G)$ in $\mfB(E^{2})$.  In general, $VN(G)$ is only a Banach algebra.  We now discuss all of this in detail.  We start with the right regular antirepresentation of $C_{c}(G)$ on $E^{2}$.  

\begin{proposition}         \label{prop:rphi}
For $F\in C_{c}(G)$, define $R_{F}:C_{c}(G)\to C_{c}(G)$ by right convolution:
$R_{F}f=f*F$.  Then $R_{F}$ extends to an element of $\mfL(E^{2})$ whose norm is $\leq \norm{F}_{I}$, and the map $F\to R_{F}$ is a $^{*}$-antirepresentation of $C_{c}(G)$, the  closure of whose image in $\mfL(E^{2})$ is canonically isomorphic to $C_{red}^{*}(G)$.
\end{proposition}
\begin{proof}
To prove that $\norm{R_{F}}\leq \norm{F}_{I}$, it is sufficient by Proposition~\ref{prop:ccgdense} to show that for $f,g\in C_{c}(G)$ and $u\in G^{0}$, that 
\begin{equation}
\mid\lan R_{F}f,g\ran(u)\mid\leq \norm{F}_{I}\norm{f}\norm{g}.       \label{eq:phifgu}
\end{equation}
This follows from (cf. \cite[p.53]{rg}):
\begin{align*}
\begin{split}
\mid\lan R_{F}f,g\ran(u)\mid & \leq 
\int\!\!\int\mid g(x)\mid\mid F(t^{-1}x)\mid\mid f(t)\mid\,d\la^{u}(t)d\la^{u}(x)  \\
& =  \int\!\!\int[\mid g(x)\mid\mid F(t^{-1}x)\mid^{1/2}]
[\mid f(t)\mid\mid F(t^{-1}x)\mid^{1/2}]\,d\la^{u}(t)d\la^{u}(x)  \\
& \leq  AB 
\end{split}
\end{align*}
where
\[ A=[\int\!\!\int\mid g(x)\mid^{2}\mid F(t^{-1}x)\mid 
\,d\la^{u}(t)d\la^{u}(x)]^{1/2},\hspace{.2in}
B=[\int\!\!\int\mid f(t)\mid^{2}\mid F(t^{-1}x)\mid\,d\la^{u}(t)d\la^{u}(x)]^{1/2}. \]
Now
\begin{align*}
A^{2}&=
\int\mid g(x)\mid^{2}\,d\la^{u}(x)\int\mid F(t^{-1}x)\mid\,d\la^{u}(t) \\
&\leq \int\mid g(x)\mid^{2}\,d\la^{u}(x)\int\mid F^{\vee}(w)\mid\,d\la^{s(x)}(w)\\
&\leq \norm{g}^{2}\norm{F^{\vee}}_{I,r}  \\
&= \norm{g}^{2}\norm{F}_{I,s}.
\end{align*}
Similarly
\[ B^{2}\leq
\norm{f}^{2}\norm{F}_{I,r}      \]
and (\ref{eq:phifgu}) follows.

Next
\begin{align*}
\lan R_{F}f,g\ran(u) &= \int\!\!\int \ov{f(t)F(t^{-1}x)}g(x)\,d\la^{u}(t)d\la^{u}(x)  \\
&= \int \ov{f(t)}\,d\la^{u}(t)\int g(x)F^{*}(x^{-1}t)\,d\la^{u}(x)     \\
&= \lan f,R_{F^{*}}g\ran(u).     
\end{align*}
So $R_{F}\in \mfL(E^{2})$, and the map $F\to R_{F}$ is a $^{*}$-antirepresentation of 
$C_{c}(G)$ into $\mfL(E^{2})$. It remains to show that for $F\in C_{c}(G)$, we have
$\norm{R_{F}}=\norm{F}_{red}$.  To prove this, for each $u$, there is an antirepresentation $R_{F}^{u}$ of $C_{c}(G)$ on $L^{2}(G^{u})$ given by right convolution by
$F$ on $C_{c}(G^{u})$.  Of course for $f\in C_{c}(G)$, 
$R_{F}f_{\mid G^{u}}=R_{F}^{u}(f_{\mid G^{u}})$.  
The map $f\to f^{\vee}$ is a linear isometry from $L^{2}(G^{u})$ onto $L^{2}(G_{u})$ that intertwines $R_{F}^{u}$ and $\pi_{u}(F^{\vee})$ ((\ref{eq:reg})).
So $\norm{R_{F}}=\sup_{u\in G^{0}}\norm{\pi_{u}(F^{\vee})}=\norm{F}_{red}$ using a well-known characterization of $\norm{F}_{red}$ (\cite{connesbook,muhlyTCU} - for more 
details, see \cite[p.108]{pat2}.)  
\end{proof}

If $\xi\in E^{2}$ and $F\in C_{c}(G)$, then the convolution formula 
\[  \xi*F(x)=\int \xi(t)F(t^{-1}x)\,d\la^{r(x)}(t)     \]
makes sense by the Cauchy-Schwartz inequality.  (The value of $\xi*F(x)$ is the same whichever representative of $\xi$ we take in the integral.)
We would expect that $\xi*F$ should be the same as $R_{F}\xi$ and be continuous on $G$, as indeed it is in the group case 
(\cite[(20.14)]{HR1}).  We now show that this is the case.

\begin{proposition}      \label{prop:cont}
Let $\xi\in E^{2}$ and $F\in C_{c}(G)$.  Then $R_{F}\xi=\xi*F$, and is a continuous function on 
$G$.  Further, if $\xi_{n}\to \xi$ in $E^{2}$, then $\xi_{n}*F\to \xi*F$ uniformly on $G$.
\end{proposition}
\begin{proof}
Let $\{f_{n}\}$ be a sequence in $C_{c}(G)$ such that $\norm{f_{n}-\xi}\to 0$.
Then given $x\in G$ and any $n\geq 1$, we have, using Proposition~\ref{prop:cont}, 
\begin{align*}
\mid (R_{F}f_{n}-\xi*F)(x)\mid & =
\mid\int (f_{n}-\xi)(t)F(t^{-1}x)\,d\la^{r(x)}(t)\mid  \\
&\leq \norm{f_{n}-\xi}^{r(x)}(\int\mid F^{\vee}(x^{-1}t)\mid^{2}\,d\la^{r(x)}(t))^{1/2}  \\
&\leq \norm{f_{n}-\xi}\norm{ F^{\vee}} \\
&\to 0
\end{align*}
independently of $x$, and $\xi*F$ is the uniform limit of a sequence of continuous functions.  
By the continuity of $R_{F}$, $R_{F}\xi=\xi*F$.  The proof of the last assertion of the proposition is similar.
\end{proof}

In the next proposition, we note that $C_{c}(G)$ is a normed algebra under the $(I,r)$-norm.  So the proposition shows that $\phi\to L_{\phi}$ is a norm decreasing homomorphism from 
$(C_{c}(G),\norm{.}_{I,r})$ into $\mfB(E^{2})$.

\begin{proposition}        \label{prop:lphi}
Let $F\in C_{c}(G)$ and $L_{F}:C_{c}(G)\to C_{c}(G)$ be the map defined by left convolution by $F$:  
$L_{F}f=F*f$.  Then 
\bi
\item[(i)] $L_{F}$ extends to a bounded linear map, also denoted by $L_{F}$, on $E^{2}$ for which $\norm{L_{F}}\leq \norm{F}_{I,r}$, and the map $F\to L_{F}$ is a norm decreasing homomorphism from $(C_{c}(G),\norm{.}_{I,r})$ into $\mfB(E^{2})$;
\item[(ii)] if $\xi,\eta\in E^{2}$ and $F\in C_{c}(G)$, then 
$F*(\xi,\eta)=(\xi,L_{F}\eta)\in B(G)$.
\ei
\end{proposition}
\begin{proof}
(i) The only non-trivial thing to be shown is that $\norm{L_{F}}\leq \norm{F}_{I,r}$.  This is equivalent to showing that for 
$f\in C_{c}(G)$ and $h^{u}\in C_{c}(G^{u})$, we have 
\begin{equation}       \label{eq:lphifh}
\mid\int_{G^{u}} \ov{(L_{F}f)(x)}h^{u}(x)\,d\la^{u}(x)\mid \leq \norm{F}_{I,r}\norm{f}
\norm{h^{u}}_{2}.
\end{equation}
To this end, 
\begin{align*}
\mid\int \ov{L_{F}f(x)}h^{u}(x)\,d\la^{u}(x)\mid &\leq
\int \left|F(t)\right|\,d\la^{u}(t)\int\left|h^{u}(x)f(t^{-1}x)\right|\,d\la^{u}(x) \\
&\leq  
\int\left|F(t)\right|\,d\la^{u}(t)\norm{h^{u}}(\int\left|f(t^{-1}x)\right|^{2}\,d\la^{u}(x))^{1/2}  \\
&\leq \int\left| F(t)\right|\norm{h^{u}}(\int\left|f(y)\right|^{2}
\,d\la^{s(t)}(y))^{1/2}\,d\la^{u}(t) \\
&\leq \norm{F}_{I,r}\norm{f}\norm{h^{u}}.
\end{align*}

(ii) If $g\in C_{0}(G)$ then 
\begin{equation}    \label{eq:FgFg}
\norm{F*g}_{\infty}\leq \norm{F}_{I,r}\norm{g}_{\infty}.
\end{equation}
It follows that the
convolution $F*(\xi,\eta)$ is defined (and continuous). Suppose first that $\eta\in C_{c}(G)$.
Then 
$(F*(\xi,\eta))(x)=\int F(t)(L_{t^{-1}x}\xi(s(x)),\eta(r(x)))\,d\la^{r(x)}(t)=\iint F(t)\ov{\xi(x^{-1}s)}
\eta(t^{-1}s)\,d\la^{r(x)}(t)d\la^{r(x)}(s)=(\xi,L_{F}\eta)$.
The same equality when $\eta\in E^{2}$ follows from (\ref{eq:FgFg}) and 
Proposition~\ref{prop:inftyxieta}.
\end{proof}

We note that if $F,f\in C_{c}(G)$ and $u\in G^{0}$, then
\begin{align*}
\lan L_{F}f,f\ran(u)&=\int \ov{F*f(x)}f(x)\,d\la^{u}(x)\\
&=\iint \ov{F(t)}\ov{f(t^{-1}x)}f(x)\,d\la^{u}(t)d\la^{u}(x)\\
&=\int \ov{F}(f,f)\,d\la^{u}.   
\end{align*}
Let $\xi\in E^{2}$ and $\{f_{n}\}$ be a sequence in $C_{c}(G)$ such that $\norm{f_{n}-\xi}\to 0$.  Taking limits in the preceding equalities with $f_{n}$ in place of $f$ then gives
\begin{equation}   \label{eq:gxixi}
\lan L_{F}\xi,\xi\ran(u)=\int \ov{F}(\xi,\xi)\,d\la^{u}.
\end{equation}
which entails in the obvious way that for $\eta\in E^{2}$,  
\begin{equation}   \label{eq:gxieta}
\lan L_{F}\xi,\eta\ran(u)=\int \ov{F}(\xi,\eta)\,d\la^{u}.
\end{equation}

\begin{proposition}           \label{prop:aident}
There exists a bounded left approximate identity $\{F_{n}\}\geq 0$ in the normed algebra
$(C_{c}(G),\norm{.}_{I,r})$, such that $L_{F_{n}}\to I$ in the strong operator topology of $\mfB(E^{2})$.
\end{proposition}
\begin{proof}
The proof is a slight modification of the proof of \cite[Proposition 1.9, p.56]{rg}.
There is a sequence $\{U_{n}\}$ of open neighborhoods of $G^{0}$ in $G$ such that each $U_{n}$ is $s$-compact
(i.e. $U_{n}\cap s^{-1}(K)$ is relatively compact for every compact subset $K$ of $G^{0}$)
and is a {\em fundamental sequence} for $G^{0}$ in the sense that
every neighborhood $V$ of $G^{0}$ in $G$ contains $U_{n}$ eventually.  There is then an increasing sequence $\{K_{n}\}$ of compact subsets of $G^{0}$ such that 
$\cup K_{n}=G^{0}$.  Using \cite[Lemma 2.12]{MRW}, there exists $g_{n}\geq 0$ in $C_{c}(U_{n})$ such that 
$\int g_{n}\,d\la^{u}=1$ for all $u\in K_{n}$. Next, there exists an open neighborhood $W_{n}$ of $K_{n}$ in $G^{0}$ such that 
$1/2< \int g_{n}\,d\la^{u}<2$ for all $u\in W_{n}$.  Let 
$h_{n}\in C_{c}(W_{n})\subset C_{0}(G^{0})$ be such that
$0\leq h_{n}\leq 1$ and $h_{n}=1$ on $K_{n}$, and set $F_{n}=g_{n}h_{n}\in C_{c}(G)$.
Then $\norm{F_{n}}_{I,r}\leq 2$. 

We now show that for $f\in C_{c}(G)$, we have $\norm{F_{n}*f-f}_{I,r}\to 0$.  Let $K$ be the support of $f$, $L$ be the (compact) closure of $U_{1}K$ in $G$, and $\eps>0$.  Then for large enough $n$, $F_{n}*f, f$ have supports inside $L$ and $\mid F_{n}*f(x)-f(x)\mid\leq \eps$ for all $x\in L$.  It follows that
$\norm{F_{n}*f - f}_{I,r}\leq \eps\sup_{u\in r(L)}\la^{u}(L\cap G^{u})$
for large enough $n$, and $\{F_{n}\}$ is a bounded left approximate identity in $C_{c}(G)$.  Similarly, $\norm{L_{F_{n}}f-f}\to 0$ in $E^{2}$ for all $f\in C_{c}(G)$.
The rest of the proposition now follows using Proposition~\ref{prop:lphi} and Proposition~\ref{prop:ccgdense}.
\end{proof}

In the preceding proposition, I do not know if there exists a bounded two-sided approximate 
identity in $C_{c}(G)$ for the $(I,r)$-norm.  It is shown in \cite[Corollary 2.11]{MRW} that there is always a two-sided (self-adjoint) approximate identity in $C_{c}(G)$ for the inductive limit topology.


\section{The Fourier algebra A(G)}

It is natural to enquire how the Fourier algebra $A(G)$ should be defined.  By analogy with the group case and also with the case of a measured groupoid (\cite{ren}), one might be inclined to take this algebra to be the closure of the span of the coefficients of $E^{2}$ in $B(G)$.  Oty (\cite[p.186]{oty}) suggests taking $A(G)$ to be the closure of $\mathcal{B}(G)\cap C_{c}(G)$ in $\mathcal{B}(G)$.  It will be convenient for our purposes to take $A(G)$ to be the closure in $B(G)$ of the subalgebra generated by the coefficients of $E^{2}$.  I do not know if the three versions of $A(G)$ coincide. 

Let $A_{cf}(G)$ be the set of coefficients $(f,g)$ of $E^{2}$ with $f,g\in C_{c}(G)$.  Note that
\begin{equation}  \label{eq:fggf}
(f,g)=\ov{g}*f^{*}.
\end{equation}
Let $A_{sp}(G)$ be the complex vector subspace of $C_{c}(G)$
spanned by $A_{cf}(G)$ and $A_{c}(G)$ be the subalgebra of $C_{c}(G)$ generated by $A_{cf}(G)$
(pointwise product).  If $V$ is an open subset of $G$, then we set
$A_{cf}(V)=A_{cf}(G)\cap C_{c}(V)$. Similarly we define $A_{sp}(V), A_{c}(V)$.
\vspace*{.2in}

\noindent
{\bf Definition}
The closure of $A_{c}(G)$ in $B(G)$ is called the {\em Fourier algebra} of 
$G$, and is denoted by $A(G)$.
\vspace*{.2in}

If $G$ is a locally compact group, then (\cite{eymard}) 
$A_{cf}(G)=A_{sp}(G)=A_{c}(G)$.  I do not know if this is true for locally compact groupoids in general.  However, when $G$ is r-discrete we have the following result.  (The discrete group version of this argument appears in \cite[Lemma VII.2.7]{davidson}.)

\begin{proposition}      \label{prop:gode}
Let $G$ be r-discrete and $\phi\in P(G)\cap C_{c}(G)$.  Then $\phi$ is a coefficient of $E^{2}$.
\end{proposition}
\begin{proof}
Let $T=R_{\phi}\in \mfL(E^{2})$.  For $g\in C_{c}(G), u\in G^{0}$, we have 
\[ \lan Tg,g\ran(u)=\ov{\iint \phi(y^{-1}x)g(y)\ov{g(x)}\,d\la^{u}(y)\,d\la^{u}(x)}\geq 0  \]
since $\phi$ is positive definite.  So $\lan T\eta,\eta\ran\geq 0$ for all $\eta\in E^{2}$, and it follows from the second part of the proof of \cite[Proposition 6.1]{paschke} that $T\geq 0$ in $C_{red}^{*}(G)\subset \mfL(E^{2})$.  Let $h\in C_{c}(G^{0})\subset C_{c}(G)$ be such that $h=1$ on 
$r(S(\phi))\cup s(S(\phi))$.  Then 
$Th(x)=\int h(t)\phi(t^{-1}x)\,d\la^{u}(t)=\phi(x)$. Further, for $x\in G$,
$(h,\phi)(x)=\int \ov{h(x^{-1}t)}\phi(t)\,d\la^{u}(t)=\phi(x)$.  Let
$\xi=T^{1/2}h\in E^{2}$.  Since $L_{F}$ commutes with every $R_{f}$, it follows from 
Proposition~\ref{prop:rphi} that it commutes with every operator $C_{red}^{*}(G)$ and hence with $T^{1/2}$.
Then for each $u\in G^{0}$ and $g\in C_{c}(G)$, we have using (\ref{eq:gxieta}), 
\begin{align*}
\int \ov{g(x)}(\xi,\xi)(x)\,d\la^{u}(x)&= \lan L_{g}T^{1/2}h,T^{1/2}h\ran(u)\\
 &=  \lan T^{1/2}L_{g}h,T^{1/2}h\ran(u)\\
&=\lan L_{g}h,Th\ran(u)\\
&=\int \ov{g(x)}(h,\phi)(x)\,d\la^{u}(x)\\
&= \int \ov{g(x)}\phi(x)\,d\la^{u}(x). 
\end{align*}
So $\phi=(\xi,\xi)$.
\end{proof}
\begin{corollary}           \label{cor:gode}
If $G$ is r-discrete then then every $\phi\in A_{c}(G)$ is the sum of two coefficients of $E^{2}$.
\end{corollary}
\begin{proof}
Let $\phi\in A_{c}(G)$. By definition, $\phi$ is a finite sum of functions that are products $\phi_{1}\ldots \phi_{n}$ where each $\phi_{i}\in A_{cf}(G)$.  By 
the construction of \cite[Proposition 1.3, (ii)$\Rightarrow$ (iii)]{rg}, there exists 
$F\in P(G\x I_{2})$ of the form (\ref{eq:bmatr})
with $\rho, \tau\in P(G)\cap C_{c}(G)$.
So $F\in C_{c}(G\x I_{2})$, and by Proposition~\ref{prop:gode}, 
$F\in A_{cf}(G\x I_{2})$.  So there exists
$\zeta\in E^{2}(G\x I_{2})$ such that $F=(\zeta,\zeta)$.  Let $\zeta_{ij}(x)=\zeta(xij)$.  Then 
$\phi=(\zeta_{21},\zeta_{11}) + (\zeta_{22},\zeta_{12})$. 
\end{proof}
\begin{corollary}  \label{cor:ideal}
If $G$ is r-discrete, then $A(G)$ is an ideal in $B(G)$.
\end{corollary}
\begin{proof}
This follows since $P(G)(A_{c}(G)\cap P(G))\subset P(G)\cap C_{c}(G)$.
\end{proof}

We now discuss another possible version $\mathcal{A}(G)$
of the Fourier algebra that coincides with $A(G)$ in the group case and relates usefully to $VN(G)$.  In fact 
$\mathcal{A}(G)$ will be a subspace of our $A(G)$ in general, and they may even be the same.  
In the group case, one way of defining $A(G)$ is to regard it as a quotient space of $L^{2}(G)\widehat{\otimes} L^{2}(G)$ (\cite{eymard,herz}).  (See also \cite[p.185]{pat1}.)  This is just the norm on $A(G)$ that comes from the identification of $A(G)$ with $VN(G)_{*}$.  This approach can be adapted, as we will see, to work for locally compact groupoids in general, with the Hilbert $D$-module $E^{2}$ replacing the $L^{2}(G)$ of the group case.

More precisely, define a map
$\theta:C_{c}(G)\x C_{c}(G)\to C_{0}(G)$ by: $\theta((f,g))=g*f^{\vee}$.  Then $\theta$ is bilinear,
and $\norm{\theta((f,g))}\leq \norm{f}\norm{g}$ by Proposition~\ref{prop:inftyxieta}.  So 
$\theta$ extends to a norm decreasing linear map, also denoted by $\theta$, from $C_{c}(G)\otimes C_{c}(G)$ (with the projective tensor product norm)
into $C_{c}(G)$.  By Proposition~\ref{prop:ccgdense}, $\theta$ extends to a norm-decreasing linear map from $E^{2}\hat{\otimes} E^{2}$ into $C_{0}(G)$.  
\vspace*{.2in}

\noindent
{\bf Definition}
The Banach space $\mathcal{A}(G)$ is defined to be the completion of the normed space
$E^{2}\hat{\otimes} E^{2}/\ker\theta$ under the quotient norm.  The norm on $\mathcal{A}(G)$ is denoted by $\norm{.}_{1}$ (this norm being like a trace class norm).
\vspace*{.2in}

We can regard $\mathcal{A}(G)$ as a linear subspace of $C_{0}(G)$.  By construction,
$A_{sp}(G)$ is a dense subspace of $\mathcal{A}(G)$.  Using the latter fact and 
Proposition~\ref{prop:inftyxieta} gives the next proposition.

\begin{proposition}        \label{prop:mcag}
$\mathcal{A}(G)$ is a subspace of $B(G)$, and for $\phi\in \mathcal{A}(G)$,  we have
\begin{equation}   \label{eq:infleq}
\norm{\phi}_{\infty}\leq \norm{\phi}\leq \norm{\phi}_{1}.
\end{equation}
\end{proposition}

The norm $\norm{.}_{1}$ on $\mathcal{A}(G)$ is given by:
\begin{equation}      \label{eq:norm2}
\norm{\phi}_{1}=\inf{\sum_{n=1}^{\infty}\norm{f_{n}}\norm{g_{n}}}
\end{equation}
the $\inf$ being taken over all expressions of the form 
$\phi=\sum_{n=1}^{\infty}g_{n}*f_{n}^{*}$ in $C_{0}(G)$ where $f_{n}, g_{n}\in C_{c}(G)$.
\vspace*{.2in}

\noindent
{\bf Note}
It would be reasonable to use a fibered projective tensor product norm in place of the projective tensor product norm in the above argument.  Indeed let $R$ be the orbit equivalence relation on $G^{0}$: so $u\sim v$ if and only if there exists an $x\in G$ such that $s(x)=u, r(x)=v$. 
Then we could consider the Banach space of functions $\phi:G\to \C$ of the form 
\begin{equation}        \label{eq:phiag}
\phi=\sum_{n=1}^{\infty} (f_{n},g_{n})   
\end{equation}
where $f_{n},g_{n}\in C_{c}(G)$ and are such that 
$M=\sup_{(u,v)\in R}(\sum_{n=1}^{\infty}\norm{f_{n}}^{u}\norm{g_{n}}^{v})<\infty$.  Further
we take $\norm{\phi}_{R}$ to be the $\inf$ of such numbers $M$.  (One uses the fact  
that for each $x\in G$, 
$\left| g*f^{\vee}(x)\right|
\leq \norm{f}^{s(x)}\norm{g}^{r(x)}$.)  I do not know if this Banach space coincides with
$\mathcal{A}(G)$.

The proof of the next proposition is left to the reader.

\begin{proposition}         \label{prop:agba}
$A(G)$ is a commutative Banach algebra, and is a subalgebra of $B(G)\cap C_{0}(G)$.  
Further, $\mathcal{A}(G)\subset A(G)$.
\end{proposition}

In our definition of the spaces $A(G), \mathcal{A}(G)$ we used the $r$-fibered Hilbert bundle $E^{2}$.  We now show that these spaces are the same if we had used the corresponding definitions using $s$ rather than $r$.  This is equivalent to saying that $A(G)=A(G(r)), 
\mathcal{A}(G)=\mathcal{A}(G(r))$ where $G(r)$ is $G$ with reversed multiplication.

\begin{proposition}    \label{prop:ggr}
$B(G)=B(G(r)), A(G)=A(G(r))$ and 
$\mathcal{A}(G)=\mathcal{A}(G(r))$.
\end{proposition}
\begin{proof}
The map $(\{H_{u}\},L)\to (\{H_{u}\},L')$, where $L'(x)=L(x^{-1})$, is a bijection from $G$-Hilbert bundles onto $G(r)$-Hilbert bundles.  It follows that the map $\phi\to \phi^{\vee}$ is an isometric isomorphism from $B(G)$ onto $B(G(r))$.  From (i) of Proposition~\ref{prop:eleag}, $B(G)=B(G(r))$.  Next, the map $f\to f^{\vee}$ is an isometry from $E^{2}(G)$ onto $E^{2}(G(r))$ and if $(f,g), (f,g)'$ denote coefficients ($f,g\in C_{c}(G)$) evaluated for $E^{2}, G$ and 
$E^{2}(G(r)), G(r)$ respectively, then $(f,g)^{\vee}=(f^{\vee},g^{\vee})'$.  It then follows that 
that $A(G)=A(G(r))$ and $\mathcal{A}(G)=\mathcal{A}(G(r))$.
\end{proof}

\begin{proposition}       \label{prop:eleag2}
\bi
\item[(i)] $\mathcal{A}(G)$ is a $D$-submodule of $C_{0}(G)$, and if $\phi\in \mathcal{A}(G)$, then $\phi^{*}\in \mathcal{A}(G)$ and $\norm{\phi}_{1}=\norm{\phi^{*}}_{1}$.  The corresponding result holds equally for $A(G)$.
\item[(ii)] If $K\in \mathcal{C}(G)$ and $U$ is an open subset of
$G$ such that $K\subset U$, then there exists $\phi\in A_{cf}(G)$ such that $\phi\in C_{c}(U)$, and $\phi(K)=\{1\}$ and $0\leq \phi\leq 1$.
\ei
\end{proposition}
\begin{proof}
(i) Use Proposition~\ref{prop:eleag}. 
(ii) (cf. \cite[Lemma 1.3]{ren}, \cite[Lemme 3.2]{eymard}) We can suppose that $\ov{U}$ is compact.  Let $L\in \mathcal{C}(G)$ be such that 
$K\subset L^{o}\subset L\subset U$.  There exists a relatively compact open subset $V$ of $G$ such that $s(L)\subset V, \ov{LV^{-1}}\subset U$ and $\ov{KV}\subset L^{o}$.  Let $f\in C_{c}(L^{o})$ be such that $0\leq f\leq 1$ and $f=1$ on $KV$.  Next let $g\in C_{c}(V)$ be such that $0\leq g\leq 1$ and $g(u)>0$ for all
$u\in s(K)$.  Let $b\in C_{c}(G^{0})$ be such that $b(u)=\la^{u}(g)$ for $u\in s(K)$, 
$b(u)\geq \la^{u}(g)$ for all $u\in G^{0}$ and $b(u)>0$ for $u\in \ov{s(U)}$.
Take $\phi=b(f*g^{\vee})\in A_{cf}(G)$.  We check that $\phi$ has the desired properties.  Obviously, $\phi\geq 0$.  Suppose that for some $x\in G$, we have $\phi(x)>0$.  Then for some $t$, $f(t), g(x^{-1}>0$, and it follows that $x\in L^{o}V^{-1}\subset 
\ov{L^{o}V^{-1}}\subset U$.  So $\phi\in C_{c}(U)$.  If $x\in K$, then any $t$ for which 
$f(t)>0, g(x^{-1}t)>0$ belongs to $KV$ so that $f(t)=1$.  (There is always such a $t$: $t=x$ will do.)  So $\phi(x)=(\int g(x^{-1}t)\,d\la^{r(x)}(t))/b(s(x))=1$.  Lastly, for any $x\in U$, 
$\phi(x)\leq (\int g(x^{-1}t)\,d\la^{r(x)}(t))/b(s(x))\leq 1$.
\end{proof}

\vspace*{.2in}

\noindent
{\bf Definition}
The set of elements $T\in \mfB(E^{2})$ such that $TR_{F}=R_{F}T$ for all $F\in C_{c}(G)$ is denoted by $VN(G)$.
\vspace*{.2in}

We now prove groupoid versions of results about $VN(G)$ proved by Eymard (\cite{eymard})
in the group case.

\begin{proposition}        \label{prop:ccont}
Let $T\in VN(G)$, $\phi\in A_{cf}(G)$ and $b\in C_{0}(G^{0})$.  Then:
\bi
\item[(i)] $T\phi$ is continuous on $G$;
\item[(ii)] $T(b\phi)=bT(\phi)$.
\ei
\end{proposition}
\begin{proof}
Write $\phi=f*g$ for $f,g\in C_{c}(G)$. Then $T\phi=TR_{g}(f)
=R_{g}(Tf)$ which is continuous by Proposition~\ref{prop:cont}.
Next, using Proposition~\ref{prop:eleag} and Proposition~\ref{prop:cont}, we get
$T(b\phi)=T(f*bg)=R_{bg}T(f)=T(f)*bg=bT(\phi)$.
\end{proof}

\begin{proposition}       \label{prop:soclos}
$VN(G)$ is the strong operator closure in $\mfB(E^{2})$ of the subalgebra 
$\mfL =\{L_{F}:F\in C_{c}(G)\}$.
\end{proposition}
\begin{proof}
Since $L_{F}\in VN(G)$ for all $F\in C_{c}(G)$ and $VN(G)$ is strong operator closed in $\mfB(E^{2})$, we have $\mfL\subset VN(G)$.  Conversely, let 
$T\in VN(G)$ and $\{F_{n}\}$ be an $(I,r)$-bounded left approximate identity for 
$C_{c}(G)$ (Proposition~\ref{prop:aident}).
Let $f_{1},\ldots ,f_{r}\in C_{c}(G)$, $N$ be a positive integer and $g\in C_{c}(G)$.  Then 
\begin{align*}
\norm{Tf_{i}-L_{g}f_{i}} & \leq \norm{T(f_{i}-L_{F_{N}}(f_{i}))} + 
\norm{(TF_{N}-g)*f_{i}} \\
&\leq \norm{T}[\max_{i}\norm{f_{i}-L_{F_{N}}(f_{i})}]+[\max_{i}\norm{R_{f_{i}}}]\norm{TF_{N}-g}.
\end{align*}
The proposition now follows using Proposition~\ref{prop:aident} and 
Proposition~\ref{prop:ccgdense} by taking $N$ large enough and $g$ close enough to $TF_{N}$ in $E^{2}$.
\end{proof}

\begin{proposition}       \label{prop:tfgu}
Let $T\in VN(G)$ and $f,g\in C_{c}(G)$.  Then 
\begin{equation}         \label{eq:tfgu}
\ov{T(f*g^{*})}=\lan Tf,g\ran,
\end{equation}
or equivalently, $T(f*g^{*})=\lan g,Tf\ran$.
\end{proposition}
\begin{proof}  
Let $\{F_{n}\}, \{U_{n}\}$ be as in the proof of of Proposition~\ref{prop:aident}.  Let
$u\in G^{0}$.  By Proposition~\ref{prop:ccont}, $T(f*g^{*})$ is continuous.
Let $\eps>0, u\in G^{0}$.  Since the sequence $\{U_{n}\}$ is fundamental, 
there exists $N$ such that for all $n\geq N$, $\mid T(f*g^{*})(x)-T(f*g^{*})(u)\mid<\eps$ for all $x\in U_{n}\cap G^{u}$.  Note also that $F_{n}\geq 0$ and $\int F_{n}\,d\la^{u}=1$. Further,
$\lan \ov{T(f*g^{*})},F_{n}\ran = \lan F_{n},T(f*g^{*})\ran = \lan R_{g}F_{n},Tf\ran
=\lan L_{F_{n}}g,Tf\ran \to \lan g,Tf\ran$ uniformly on $G^{0}$.  Now 
\begin{align*}
\mid \lan \ov{T(f*g^{*})},F_{n}\ran(u) - T(f*g^{*})(u)\mid &
\leq \int \mid T(f*g^{*})(x)-T(f*g^{*})(u)\mid F_{n}(x)\,d\la^{u}(x)  \\
&\leq \sup_{x\in U_{n}\cap G^{u}}\mid T(f*g^{*})(x)-T(f*g^{*})(u)\mid \to 0.
\end{align*}
So $\ov{T(f*g^{*})(u)}=\lan Tf,g\ran(u)$.
\end{proof}   

\begin{proposition}        \label{prop:tphix}
Let $T\in VN(G)$.  Then $T$ determines an element, also denoted by $T$, in $\mfB(\mathcal{A}(G))$, and its norm $\norm{T}_{1}$ in $\mfB(\mathcal{A}(G))$ is $\leq \norm{T}$.
\end{proposition}
\begin{proof}
Let $\phi\in A_{sp}(G)$.  Then $\phi\in E^{2}$ so that $T\phi\in E^{2}$.  If $\phi=g*f^{*}$ 
($f,g\in C_{c}(G)$) then by Proposition~\ref{prop:cont},
$T\phi=Tg*f^{*}\in \mathcal{A}(G)$.
Further, $\norm{T\phi}_{1}\leq \norm{T}\norm{g}\norm{f}$, and $\norm{T}_{1}\leq \norm{T}$.
\end{proof}
\begin{corollary}     \label{cor:tfgu}
If $\phi\in \mathcal{A}(G)$, then $\norm{T\phi}_{\infty}\leq \norm{T}\norm{\phi}_{1}$ for all $x\in G$.
\end{corollary}
\begin{proof}
Use (\ref{eq:infleq}).
\end{proof}

\noindent
{\bf Example}\\
\noindent
Here is a very simple example to show (among other things) that 
in the situation of the Proposition~\ref{prop:lphi}, we 
do not usually obtain that $L_{F}\in \mfL(E^{2})$.  Let $G=X\x X$ be a trivial groupoid  with measure $\mu$ on $X$.  We can (and indeed will) take $G$ to be $G_{n}=\{1,2,\ldots ,n\}\x \{1,2,\ldots ,n\}$ and $\mu$ counting measure on $\{1,2,\ldots ,n\}$. 
Note that for $(x,y)\in G$, $r(x,y)=x$ and $s(x,y)=y$.  

For $F\in C_{c}(G)$, we have 
\[R_{F}f(x,y)=(f*F)(x,y)
=\int f(x,t)F((t,x)(x,y))d\mu(y)=\int f(x,t)F(t,y)\,d\mu(t).  \]
Also, 
$L_{F}f(x,y)=\int F(x,t)f(t,y)\,d\mu(t)$.  

Now let $G=G_{n}$.  We will calculate $C_{red}^{*}(G)$, $VN(G)$, $\mfL(E^{2})$, $VN(G)\cap C_{red}^{*}(G)$, $\mathcal{A}(G)$ and   
$A(G)$.  We identify $G^{0}$ with $\{1,2,\ldots ,n\}$ and $G^{i}=\{i\}\x G^{0}$ with $G^{0}$.
We also identify $E^{2}$ with $M_{n}$, where for $f\in M_{n}$, the function $f^{i}$ $(i\in G^{0})$ on 
$G^{0}$ is given by: $f^{i}(j)=f_{ij}$.  The inner product on $E^{2}$ is given by:
$\lan A,B\ran(i)=\sum_{j}\ov{A_{ij}}B_{ij}$.  Trivially,
$C_{red}^{*}(G)=C_{c}(G)$ is just $M_{n}$ multiplying itself on the right, and the adjoint of
$A\in M_{n}$ is the usual adjoint $A^{*}$.  The elements $T$ of $C_{red}^{*}(G)$ are thus those for which there is a matrix $\psi$ such that 
\begin{equation}
T(e_{ij})=\sum_{l} \psi_{lj}e_{il}.     \label{eq:credpsi}
\end{equation}
Similarly, $VN(G)$ is just $M_{n}$ multiplying itself on the left.  

Now let $T\in \mfB(E^{2})$.  Write $T(e_{ij})=\sum \al_{ijkl}e_{kl}$.
Suppose that $T$ is a module map.  Then for all $b\in C(X)$, we have 
\[  \sum\al_{ijkl}b(i)e_{kl}=T(e_{ij}b)=(Te_{ij})b=\sum b(k)\al_{ijkl}e_{kl}.     \]
It follows that $T(e_{ij})=\sum \al_{ijil}e_{il}$, and it is easily checked that the latter is a necessary and sufficient condition for $T\in \mfB(E^{2})$ to belong to $\mfL(E^{2})$, with 
$T^{*}$ given by: $T^{*}(e_{il})=\sum \ov{\al_{ijil}}e_{ij}$.  So the dimensions of 
$\mfB(E^{2}), \mfL(E^{2}), VN(G)$, and $C_{red}^{*}(G)$ are respectively $n^{4}, n^{3}, n^{2}$
and $n^{2}$.  We now show that $VN(G)$ and $C_{red}^{*}(G)$ intersect in the multiples of the identity, so that all four spaces are different when $n>1$.

For $T$ to belong to $C_{red}^{*}(G)\cap VN(G)$, we require first that $T=L_{\phi}$ for some $\phi\in M_{n}$
so that $T(e_{ij})=\sum \phi_{ki}e_{kj}$ (and in the notation of the preceding paragraph, $a_{ijkl}=\phi_{ki}$ if $l=j$ and is $0$ otherwise).  For $T=L_{\phi}$ to belong to 
$\mfL(E^{2})$, we require by the preceding paragraph that $\phi_{ki}=0$ when $k\ne i$, so that  for all $i,j$,
$T(e_{ij})=\phi_{ii}e_{ij}$.  Comparing this with (\ref{eq:credpsi}) gives that $T$ is a multiple of the identity, so that $VN(G)\cap C_{red}^{*}(G)=\C 1$.

As vector spaces, $\mathcal{A}(G)=A(G)=B(G)=M_{n}$.  Indeed, as algebras, all four algebras are just $M_{n}$ under the Schur product (using the functions $g*f$ and pointwise multiplication).  A result of Paulsen (\cite[p.31]{paul}) shows that
$\norm{\phi}_{cb}\leq \norm{\phi}_{C^{*}(G)}$ for all $\phi\in A(G)$.  Another result of Paulsen
(\cite[p.112]{paul}) can be used to show that $\norm{.}_{cb}=\norm{.}$ on $B(G)$.
Indeed, the result is that $\norm{A}_{cb}\leq 1$ if and only if there exist $f,g\in E^{2}$ with $\norm{f}, \norm{g}\leq 1$ and $A_{ij}=(f_{j},g_{i})=(f,g)(i,j)$.  It follows that
$\norm{.}\leq \norm{.}_{cb}$ on $A(G)$, and equality follows from Corollary~\ref{cor:norms}.  It also follows that $\norm{.}=\norm{.}_{1}$.

\vspace*{.2in}
An important fact used by Eymard in his study of $A(G)$ in the group case is that 
$VN(G)$ is identifiable in the natural way with $A(G)^{*}$ (i.e. $A(G)$ is the predual of the von Neumann algebra $VN(G)$).  We need the groupoid version of this result.  Of course in the groupoid case, $VN(G)$ is not a von Neumann algebra, but despite that, we will show that there is a suitable version of this identification for groupoids. 

Let $\mfB_{D}(\mathcal{A}(G),D)$ be the Banach space of continuous, linear, right $D$-module maps from $\mathcal{A}(G)$ into $D$.  The space $\mfB_{D}(\mathcal{A}(G),D)$ is a left $D$-module with dual action: $b\al(\phi)=\al(\phi b)$.  
We think of $\mfB_{D}(\mathcal{A}(G),D)$ as the ``dual'' of $\mathcal{A}(G)$, and write it as $\mathcal{A}(G)'$. Of course, $\mathcal{A}(G)'$ is very different from the {\em Banach} space dual space $\mathcal{A}(G)^{*}$ in general, but the two do coincide in the group case.  

For $f\in C_{c}(G), \al\in \mfB_{D}(\mathcal{A}(G),D)$, let $f\al:C_{c}(G)\to D$ be given by: $f\al(g)=\al(g*f)$.  Then $f\al$ is linear, and  
$\norm{f\al(g)}\leq \norm{\al}\norm{g*f}_{1}\leq \norm{\al}\norm{g}\norm{f^{*}}$.  So $f\al$ extends to a bounded linear map, also denoted by $f\al$, from $E^{2}$ into $D$.  
Also each map $f\al$ is a (right) module map.  Indeed, using (\ref{eq:baction}) and the fact 
that $\al$ is a module map, we have 
$f\al(gb)=\al(gb*f)=\al((g*f)b)=\al(g*f)b=(f\al(g))b$.  Further, for fixed $f$, the map
$\al\to f\al$ is bounded and linear from $\mathcal{A}(G)'$ into $\mathcal{B}_{D}(E^{2},D)$.

Let $\mathcal{A}_{K}^{r}(G)'$ be the set of $\al\in \mfB_{D}(\mathcal{A}(G),D)$ such that $f\al\in \mfK(E^{2},D)$ for all $f\in C_{c}(G)$.  It follows from the continuity of each of the maps $\al\to f\al$ and the closedness of $\mfK(E^{2},D)$ in $\mfB_{D}(E^{2},D)$ that $\mathcal{A}_{K}^{r}(G)'$ is a closed subspace of $\mfB_{D}(\mathcal{A}(G),D)$.  
Further, $\mathcal{A}_{K}^{r}(G)'$ is a left invariant subspace of $\mfB_{D}(\mathcal{A}(G),D)$. 
This follows since $f(b\al)=(fb)\al\in \mfK(E^{2},D)$ whenever $\al\in \mathcal{A}_{K}^{r}(G)'$.
If $G$ is a locally compact group, then $\mathcal{A}_{K}^{r}(G)'=\mathcal{A}(G)^{*}$.  

If $G=G_{n}$, then $\mfB_{D}(\mathcal{A}(G),D)=\mathcal{A}_{K}^{r}(G)'$.  Indeed, 
let $\al\in \mfB_{\C^{n}}(M_{n},\C^{n})$.  Since $\al$ is a $\C^{n}$-module map, we have 
$\al(e_{ij})=\la_{ij}e_{i}$ for some $\la_{ij}\in \C$.  Then with $f=B\in M_{n}$, we have 
$B\al(e_{ij})=\ov{\eta_{ij}}e_{i}$, where $\eta_{ij}=\ov{\sum_{k}\la_{ik}B_{jk}}$.  Then 
$B\al=((1,1,\ldots ,1),\eta)\in \mfK(E^{2},D)$.

We can define $\mathcal{A}_{K}^{l}(G)'$ to be $\mathcal{A}_{K}^{r}(G(r))'$.

The following theorem  shows that $VN(G)$ identifies naturally with $\mathcal{A}_{K}^{r}(G)'$ as a Banach space, and
generalizes \cite[Th\'{e}or\`{e}me (3.10)]{eymard}. 


\begin{theorem}      \label{th:agdual}
For each $T\in VN(G)$, there exists a unique element $\al_{T}\in \mathcal{A}_{K}^{r}(G)'$ defined by:
\begin{equation}      \label{eq:altphi}
\al_{T}(\phi)=\ov{T(\phi^{*})}_{\mid G^{0}}.\hspace{.2in}(\phi\in \mathcal{A}(G))
\end{equation}
Further
\begin{equation}
\al_{T}(g*f^{*})(u)=\lan Tf,g\ran(u) = \ov{T(f*g^{*})(u)}.  \label{eq:alT}
\end{equation}
Lastly, the map $T\to \al_{T}$ is a linear isometry from $VN(G)$ onto $\mathcal{A}_{K}^{r}(G)'$.
\end{theorem}
\begin{proof}
Let $\al\in \mathcal{A}_{K}^{r}(G)'$.  Let $f\in C_{c}(G), u\in G^{0}$.  Since $f^{*}\al\in
\mfK(E^{2},C_{0}(G^{0}))$, there exists a unique $F_{f}\in E^{2}$ such that 
$\al(g*f^{*})=\lan F_{f},g\ran$.  (This is the Riesz-Fr\'{e}chet theorem for Hilbert
modules (\cite[p.13]{lance}).)
Define a linear operator $T$ on $E^{2}$ by setting $Tf=F_{f}$.    For any $u$, we can find 
$g\in C_{c}(G)$ such that $((Tf)^{u},g^{u})$ is close to $\norm{(Tf)^{u}}$ and both $\norm{g^{u}}$ and $\norm{g}$ close to $1$.  
Since $\norm{\lan Tf,g\ran} \leq \norm{\al}\norm{g*f^{*}}_{1}\leq\norm{\al}\norm{g}\norm{f}$,  we obtain that $T$ is bounded with
$\norm{T}\leq \norm{\al}$.  For the reverse inequality, let $\phi\in \mathcal{A}(G)$.  Suppose first that  
$\phi=g*f^{*}$ for some $f,g\in C_{c}(G)$.  Then 
$\norm{\al(\phi)}=\norm{\lan Tf,g\ran}\leq \norm{T}\norm{f}\norm{g}$.  It follows that for 
general $\phi\in \mathcal{A}(G)$, we have 
$\norm{\al(\phi)}\leq\norm{T}\norm{\phi}$.  So $\norm{\al}=\norm{T}$.  Next we show that $T\in VN(G)$. Indeed, for $f_{1},f_{2}\in C_{c}(G)$, we have 
$\lan TR_{f_{2}}f_{1},g\ran=\al(g*(f_{1}*f_{2})^{*})=
\al(g*f_{2}^{*}*f_{1}^{*})=\lan Tf_{1},(R_{f_{2}})^{*}g\ran
=\lan R_{f_{2}}Tf_{1},g\ran$, so that $T\in VN(G)$.  (\ref{eq:altphi}) follows from Proposition~\ref{prop:tfgu}.  

Conversely, let $T\in VN(G)$.  For $\phi\in A_{sp}(G)$, define 
$\al(\phi)=\ov{T(\phi^{*})}_{\mid G^{0}}$.  Then using Corollary~\ref{cor:tfgu} and
(ii) of Proposition~\ref{prop:eleag2},
$\al$ extends to a bounded linear map on 
$\mathcal{A}(G)$, and is a right module map since
$\al(\phi b)(u)=\ov{T(\ov{b}\phi^{*})}=\al(\phi)b$ by Proposition~\ref{prop:ccont}.  
Since $f\al(g)=\lan T(f^{*}),g\ran$, we have 
that $\al\in \mathcal{A}_{K}^{r}(G)'$.  Trivially, $\al=\al_{T}$. 
\end{proof}   
\noindent
{\bf Example}\\
Let $F\in C_{c}(G)$, $T=L_{F}$ and $\al=\al_{T}\in \mathcal{A}_{K}^{r}(G)'$.  From (\ref{eq:alT}), for 
$g\in C_{c}(G)$, we have 
\[  \al(g*f^{*})(u)=\ov{L_{F}(f*g^{*})(u)}=\ov{\int F(t)(f*g^{*})(t^{-1})\,d\la^{u}(t)},  \]   
from which it follows that $\al(\phi)=\lan F,\phi\ran$.
\vspace*{.2in}


\section{Duality for $A(G)$}

In this section, we prove a groupoid version of Eymard's duality theorem for groups. Eymard's duality result says that the character space (i.e. the space of non-zero multiplicative linear functionals) on $A(G)$ is just $G$ itself.  I do not know if the character space of the commutative Banach algebra $A(G)$ ($G$ a groupoid) can be identified with $G$ as in the group case.  Instead, we replace scalar-valued homomorphisms by $D$-valued module homomorphims.  We will obtain a duality theorem for the groupoid case which coincides with Eymard's duality theorem in the group case.  We will deal initially with a more general situation than is strictly required for our main theorem since the former may prove useful for a more general duality theorem.   

For each $u\in G^{0}$, adjoin a point $\infty^{u}$ to $G^{u}$, and let $H^{u}=G^{u}\cup \{\infty^{u}\}$, $H=\cup_{u\in G^{0}}H^{u}$.  Extend the range map $r$ to $H$ by defining $r(\infty^{u})=u$.  

We give $H$ a locally compact Hausdorff topology as follows.   (Each subspace $H^{u}$ in the relative topology will turn out to be the one-point compactification of $G^{u}$.)  Let  $\mathcal{B}$ be the family of sets that are either of the form $U$ or of the form $V$, where $U$ is any open subset of $G$, and $V$ is of the form $r^{-1}(W)\setminus C\subset H$ where $W$ is any open subset of $G^{_0}$ and $C$ is a compact subset of $G$.  The proof of the following proposition is left to the reader.

\begin{proposition}        \label{prop:H}
The family $\mathcal{B}$ is a basis for a locally compact, second countable Hausdorff topology on $H$, and $r:H\to G^{0}$ is an open map.  Next, the relative topology inherited by $G$ from $H$ is the original topology of $G$, and $G$ is an open subset of $H$.  Further, for any $u\in G^{0}$ and any sequence $\{x_{n}\}$ in $H$, we have $x_{n}\to \infty^{u}$ if and only if 
for any compact subset $C$ of $G$, the sequence
$\{x_{n}\}$ is in $H\setminus C$ eventually, and $r(x_{n})\to u$.  Lastly, the map $u\to \infty^{u}$ is continuous, and the relative topology on $H^{u}$ is that of the one point compactification of $G^{u}$.
\end{proposition}

Let $\Ga^{r}$ be the set of continuous sections of $(G,r)$, i.e. the set of continuous functions 
$\ga:G^{0}\to G$ such that $\ga(u)\in G^{u}$ for all $u\in G^{0}$.  Similarly, let $\De^{r}$ be the set of continuous sections of $(H,r)$.  Of course, $\Ga^{r}\subset \De^{r}$.  Since
$\phi\in A(G)\subset C_{0}(G)$, we can regard $\phi\in C_{0}(H)$ by defining 
$\phi(\infty^{u})=0$ for all $u\in G^{0}$.  For
$\ga\in \De^{r}$, define a map $\al_{\ga}:A(G)\to C_{0}(G^{0})$ by:
\[      \al_{\ga}(\phi)=\phi\circ \ga.                      \]

An element $\al\in \mfB_{D}(A(G),D)$ is said to be {\em multiplicative} if $\al(\phi\psi)=\al(\phi)\al(\psi)$ for all $\phi,\psi\in A(G)$.  The set of multiplicative elements of \\
$\mfB_{D}(A(G),D)$ is denoted by $\Phi^{r}_{A(G)}$.  When $G$ is a group, then
$\Phi^{r}_{A(G)}$ is just the set of multiplicative linear functionals on $A(G)$.
	
\begin{proposition}              \label{prop:alga}
The map $\al_{\ga}$ belongs to $\Phi^{r}_{A(G)}$ for all $\ga\in \De^{r}$.  
\end{proposition}
\begin{proof}
Let $\phi\in A(G)$. Then (Proposition~\ref{prop:posu})
\[ \norm{\al_{\ga}(\phi)}=\sup_{u\in G^{0}}\mid\phi(\ga(u))\mid\leq \norm{\phi}.   \]
Further, $\al_{\ga}$ is a module map since
\[ \al_{\ga}(\phi b)(u)=\phi(\ga(u))b(r(\ga(u)))=(\al_{\ga}(\phi)b)(u). \]
Next, it is trivial that $\al_{\ga}$ is multiplicative on $A(G)$.  So  $\al_{\ga}$ belongs to $\Phi^{r}_{A(G)}$.
\end{proof}

I do not know if every element of $\Phi^{r}_{A(G)}$ is of the form $\al_{\ga}$ for some
$\ga\in \De^{r}$.  A tentative conjecture is that the answer is {\em yes} and that  
$\De^{r}$ with the topology of pointwise convergence corresponds to $\Phi^{r}_{A(G)}$ 
with the pointwise topology on $A(G)\x G^{0}$.  In the group case this is effectively
Eymard's theorem, except that we are allowing the 0-linear functional in the character space of $A(G)$.  (This functional corresponds to $\ga(e)=\infty^{e}$ where $e$ is the unit of the group $G$.)
In our present situation, it is reasonable to allow elements of $\Phi^{r}_{A(G)}$ to vanish for some $u$'s, i.e. to allow the existence of $u$'s for which $\al(\phi)(u)=0$ for all $\phi\in A(G)$.  Adding on the points $\infty^{u}$ allows one to incorporate this within the section viewpoint.  

Instead, our main theorem is also a generalization of Eymard's theorem in which $\De^{r}$ is replaced by elements of $\Ga^{r}$ that have a certain symmetry with respect to the $s$-map, 
and indeed correspond to bisections of $G$.  For the present, we show that every $\al\in \Phi^{r}(A(G))$ is associated with (at least) a partially defined continuous section on $G^{0}$.  

So let $\al\in \Phi^{r}_{A(G)}$.  Define
$A$ to be the set of $x\in G$ such that for every neighborhood $V$ of $x$, there exists $\phi\in A_{cf}(V)$ such that $\al(\phi)(r(x))\ne 0$.
Define $N$ to be the set of
$x\in G$ for which there exists an open neighborhood $U$ of $x$ such that for all
$\phi\in A_{cf}(U)$, we have $\al(\phi)=0$.
Trivially, $N$ is an open subset of $G$.   
Write $A^{0}=r(A)$ and $B^{0}=G^{0}\setminus A^{0}$.  

\begin{proposition}       \label{prop:AN}
\bi
\item[(i)] For every $u\in A^{0}$, the set $A\cap G^{u}$ is a singleton $\{x^{u}\}$.
\item[(ii)] For every $u\in A^{0}$, we have $G^{u}\setminus \{x^{u}\}\subset N$.
\item[(iii)] If $u\in B^{0}$ and $\phi\in A_{cf}(G)$, then $\al(\phi)(u)=0$.
\item[(iv)] If $u\in G^{0}$, $\phi\in A_{cf}(G)$ and $S(\phi)\cap G^{u}\subset N$, then  $\al(\phi)=0$ in a neighborhood of $u$ in $G^{0}$.
\item[(v)] $A^{0}$ is an open subset of $G^{0}$.
\item[(vi)] The map $u\to x^{u}$ is continuous from $A^{0}$ into $G$.
\item[(vii)] Let $\phi\in A_{cf}(G)$ be such that 
$S(\phi)\subset G\setminus \{x^{u}: u\in A^{0}\}$.  Then $\al(\phi)=0$.
\ei
\end{proposition}
\begin{proof}
(i) and (ii). Let $x\in A$ and $u=r(x)$.  Let $y\in G^{u}$ with $y\ne x$.  Let $W_{x}, W_{y}$ be disjoint open neighborhoods of $x,y$ in $G$.  Let $\phi\in A_{cf}(W_{x})$ be such that $\al(\phi)(r(x))\ne 0$.  Let $U_{x}=\{z\in W_{x}: \al(\phi)(r(z))\ne 0\}$.  Since $\al(\phi)$ is continuous, we have that $U_{x}$ is an open neighborhood of $x$.  Let 
$U_{y}=r^{-1}(r(U_{x}))\cap W_{y}$.  Since $r$ is open and continuous and $r(y)=r(x)\in r(U_{x})$, it follows that $U_{y}$ is an open neighborhood of $y$, and $r(U_{y})\subset  
r(U_{x})$.  Let $\psi\in A_{cf}(U_{y})$.  We will show that $\al(\psi)=0$ so that $y\in N$.
Indeed, since $\phi\psi=0$ (since $W_{x}\cap U_{y}\subset
W_{x}\cap W_{y}=\emp$), we have $0=\al(\phi\psi)=\al(\phi)\al(\psi)$.  Since $\al(\phi)(u)\ne 0$ on $r(U_{y})$, it follows that $\al(\psi)=0$ on $r(U_{y})$.  It remains to show that $\al(\psi)$ vanishes outside $r(U_{y})$.  To this end, $r(S(\psi))$ is a compact subset of $r(U_{y})$.
Let $b\in C_{c}(G^{0})$ be such that $b=1$ on $r(S(\psi))$ and $0$ outside $r(U_{y})$.  Then 
$b\circ r=1$ on $S(\psi)$, so that 
\[  \al(\psi)=\al(\psi b)=\al(\psi)b=0     \]
outside $r(U_{y})$.  (i) and (ii) now follow since $A\cap N=\emp$.  

(iii) Let $u\in B^{0}$ and $\phi\in A_{cf}(G)$.  Let $C^{u}=S(\phi)\cap G^{u}$.  Since $G^{u}\cap A=\emp$, we can cover $C^{u}$ by a finite number of relatively compact, open sets $U_{1},\ldots ,U_{n}$, $U_{i}\cap G^{u}\ne \emp$ for each $i$, and such that $\al(\psi)(u)=0$ for all $\psi\in A_{cf}(U_{i})$, $1\leq i\leq n$. Let $U=\cup_{i=1}^{n}U_{i}$.  There exists an open neighborhood $W$ of $u$ in $G^{0}$ and a function $b\in C_{c}(W)$ with $b(u)=1$ such that  
$\phi'=\phi b\in A_{cf}(U)$.  Since $\al$ is a module map and $b(u)=1$, we have
$\al(\phi')(u)=\al(\phi)(u)$.  Using a partition of unity argument (\cite[p.7]{helg}) we can write $\phi'=\sum_{j=1}^{m}\phi'_{j}$ where $\phi'_{j}\in A_{cf}(U_{i})$.  Then 
$\al(\phi')(u)=\sum_{j=1}^{m}\al(\phi'_{j})(u)=0$, and (iii) is proved. (iv) is proved in the same way as (iii). 

(v), (vi) and (vii). Let $u_{0}\in A^{0}$.  Let $U$ be an open neighborhood of $x^{u_{0}}$ in $G$.  Since $x^{u_{0}}\in A$, there exists $\phi\in A_{cf}(U)$ such that $\al(\phi)(u_{0})\ne 0$.  By continuity, there exists an open subset $W$ of $r(U)$ such that $u_{0}\in W$ and
$\al(\phi)(u)\ne 0$ for all $u\in W$.  Let $u\in W$.
By (iii), $u$ does not belong to $B^{0}$ and so $W\subset A^{0}$.  (v) now follows.  By (iv), $S(\phi)\cap G^{u}$ is not contained in $N$.  From (ii), $x^{u}\in S(\phi)\subset U$.  (vi) now follows since the inverse image of $U$ under the map $u\to x^{u}$ contains an open neighborhood of $u_{0}$.  (vii) follows from (ii), (iii) and (iv).
\end{proof}

Let $\al,\{x^{u}\}$ be as in Proposition~\ref{prop:AN}.  Define $\ga:G^{0}\to H$ by setting $\ga(u)=x^{u}$ if $u\in A^{0}$ and $=\infty^{u}$ otherwise.  We note that $\ga\in \Ga^{r}$ if and only if $A^{0}=G^{0}$.

\begin{proposition}      \label{prop:ga}
The map $\ga\in \De^{r}$ if and only if $G=A\cup N$.
\end{proposition}
\begin{proof}
Suppose that $\ga\in \De^{r}$. Suppose that $G\ne (A\cup N)$.  Let $x_{0}\in G\setminus (A\cup N)$ and $u_{0}=r(x_{0})$.  By Proposition~\ref{prop:AN}, (ii), (iii), $u_{0}\in B^{0}$ and for all
$\phi\in A_{cf}(G)$, $\al(\phi)(r(x_{0}))=0$.  Let $U$ be a neighborhood of $x_{0}$.  Since
$x_{0}\notin A\cup N$, there exists a $\psi\in A_{cf}(U)$ such that $\al(\psi)\ne 0$. Using (iii) and (ii) of Proposition~\ref{prop:AN}, $x^{v}\in U$ for some $v\in r(U)\cap A^{0}$.
So there exists a sequence $\{v_{n}\}$ in $A^{0}$ such that $v_{n}\to u_{0}$ and $\ga(v_{n})=x^{v_{n}}\to x_{0}$.  Since $\ga$ is continuous, $x_{0}=\ga(u_{0})$ and $x_{0}\in A$.  This is a contradiction.  So $G=A\cup N$.

Conversely, suppose that $G=A\cup N$ and let $u_{n}\to u$ in $G^{0}$.  If $u\in A^{0}$, then
by (v) and (vi) of Proposition~\ref{prop:AN}, $u_{n}\in A^{0}$ eventually, and $\ga(u_{n})\to \ga(u)$.  Suppose then that $u\in B^{0}$.  Let $C$ be a compact subset of $G$.  Suppose that $\ga(u_{n})\in C$ eventually.  (So $u_{n}\in A^{0}$ eventually.)
We can suppose that $\ga(u_{n})\to x$ for some $x\in C$.  Since $r(x)=u\notin A^{0}$ and $G=A\cup N$, it follows that $x\in N$. So $\ga(u_{n})\in N$ eventually (since $N$ is open) and we contradict $\ga(u_{n})\in A$.  So $\ga(x_{n})\notin C$ eventually, and $r(\ga(u_{n}))=u_{n}\to u$.
From Proposition~\ref{prop:H}, $\ga(x_{n})\to \infty^{u}=\ga(u)$. 
\end{proof}   

For the rest of this paper, we leave $H$ behind and instead consider only maps $\al\in \Phi^{r}_{A(G)}$ for which 
$A^{0}=G^{0}$.  Then $\ga$ is an $r$-section of $G$.  We now want $\ga$ to determine also an $s$-section of $G$.  Two reasons for this are as follows.  Firstly, we would like, as in the group case, to have a multiplicative structure on the set of such $\ga$'s and the product of two $r$-sections is not usually an $r$-section.  Another indication that we need to consider $s$-sections is that $r$-sections give right module maps, and a reasonable ``dual'' for $A(G)$ should not have a ``bias'' for right over left.  To deal with these, 
we want to pair together $\al$ as above with a corresponding {\em left} module map $\bt$.  The pair $(\al,\bt)$ can usefully be thought of in terms similar to that for multipliers, where a (two-sided) multiplier is a pairing of a left and right multiplier. 

\vspace*{.2in}
\noindent
{\bf Definition}
A {\em multiplicative module map on $A(G)$} is a pair of maps 
$\al,\bt\in \mfB(A(G),D)$ such that:
\bi
\item[(i)] $\al\in \mfB_{D}(A(G),D)$ and $\bt\in _{D}\mfB(A(G),D)$;
\item[(ii)] $\al(A(G))(u)\ne 0$ for all $u\in G^{0}$;
\item[(iii)] there exists a homeomorphism $J$ of $G^{0}$ such that for all $\phi\in A(G)$,
\begin{equation}
     \bt(\phi)\circ J=\al(\phi);       \label{eq:btalt}
\end{equation}
\item[(iv)] $\al_{\mid \mathcal{A}(G)}\in \mathcal{A}^{r}_{K}(G)'$;
\item[(v)] $\al$ is multiplicative on $A(G)$.
\ei
The set of multiplicative module maps on $A(G)$ is denoted by $\Phi_{A(G)}$.
\vspace*{.2in}

In the above definition, (ii), (iv) and (v) involve only $\al$.  But in fact using 
(i) and (iii), the corresponding properties for $\bt$ in (ii),(iv) and (v) follow.  
(In (iv), $\bt\in G(r)$.)  
So the roles of $\al, \bt$ in a multiplicative module map on $A(G)$ are symmetrical.

Let $(\al,\bt)\in \Phi_{A(G)}$. The condition that $\al(A(G))(u)\ne 0$ in (ii) is equivalent to the section $\ga$ associated with $\al$ in Proposition~\ref{prop:AN} being a global section of $(G,r)$, i.e. $A^{0}=G^{0}$.  Similarly, $\bt$ determines a section $\de$ of 
$(G,s)$.  We now show that the pair $(\ga,\de)$ determine a {\em bisection}, i.e. a subset $B$ of $G$ such that both $r:B\to G^{0}$ and $s:B\to G^{0}$ are homeomorphisms. 

\begin{proposition}         \label{prop:bisec}
There exists a bisection $B$ such that $\ga=(r_{\mid B})^{-1}$ and $\de=(s_{\mid B})^{-1}$.
\end{proposition}
\begin{proof}
Let $u\in G^{0}$ and $x=\ga(u)$.  By the definition of $\ga$, for every open neighborhood   
$V$ of $x$, there exists $\phi\in A_{cf}(V)$ such that $\al(\phi)(u)\ne 0$. Then from (iii),
$\bt(\phi)(Ju)=\al(\phi)(u)\ne 0$.  Suppose that $Ju\ne s(x)$.  By contracting $V$, we can suppose that $S(\phi)\cap G_{J(u)}=\emptyset$.  Applying (iv) of Proposition~\ref{prop:AN} to $G(r)$, it follows that $\bt(\phi)(Ju)=0$.  This is a contradiction.  So $Ju=s(x)$.  By the definition of $\de$, we have $\de(s(x))=x$.  Take $B$ to be the range of $\ga$ ($=$ the range of 
$\de$).
\end{proof}

In the above, we have $Ju=s(\ga(u))=(\de^{-1}\circ \ga)(u)$, and $\bt$ is uniquely determined by $\al$, i.e. if there is a $\bt$ such that $(\al,\bt)\in \Phi_{A(G)}$, then $\bt$ is unique. 
Let $\Ga$ be the set of bisections of $G$.  Each $B\in \Ga$ is identifiable with a pair $(\ga,\de)$ of sections of $(G,r)$ and $(G,s)$ respectively, and conversely.  
Further, $\Ga$ is a group under setwise products and inversion.  See e.g. \cite[Proposition 2.2.4]{pat2} for the r-discrete case.  To check that $BC\in \Ga$ if $B,C\in \Ga$ suppose that $B,C$ are associated with the pairs of sections $(\ga,\de)$ and $(\ga',\de')$ respectively.  Then the pair of sections associated with $BC$ are $u\to \ga(u)\ga'(s(\ga(u)))$ and $u\to \de(r(\de'(u)))\de'(u)$ respectively.  

The main theorem of this paper shows that under certain conditions, $\Phi_{A(G)}=\Ga$.  Let us make this more precise.  Each $a\in \Ga$ determines sections $\ga:u\to a^{u}$, $\de:u\to a_{u}$ as above.  These in turn determine bounded linear maps $\al^{a},\bt_{a}$ from $A(G)$ into $C_{0}(G^{0})$ by setting:
\begin{equation}           \label{eq:albta}
\al^{a}(\phi)(u)=\phi(a^{u}),     \bt_{a}(\phi)(u)=\phi(a_{u}).       
\end{equation}
\begin{proposition}            \label{prop:alabta}
The pair $(\al^{a},\bt_{a})$ belongs to $\Phi_{A(G)}$.
\end{proposition}
\begin{proof}
Firstly,  $\al^{a}$ is a bounded linear operator from $A(G)$ into $D$ using Proposition~\ref{prop:posu}. It is also a module map since for $b\in D$, we have 
$\al(\phi b)(u)=(\phi b)(a^{u})=\al(\phi)(u)b(u)=(\al(\phi)b)(u)$.  So 
$\al^{a}\in \mfB_{D}(A(G),D)$.  Next, if $f,g\in C_{c}(G)$, we have 
$f\al^{a}(g)(u)=(g*f)(a^{u})=\int g(t)f(t^{-1}a^{u})\,d\la^{u}(t)=\theta_{h,k}(g)(u)$ where 
$k(t)=\ov{f(t^{-1}a^{r(t)})}\in C_{c}(G)$ and $h\in C_{c}(G^{0})$ is such that 
$h(u)=1$ on $(\ga^{-1}\circ \de)(s(S(f)))$.  So $\al^{a}\in \mathcal{A}_{K}^{r}(G)'$. 
Similarly, $\bt_{a}\in \mathcal{A}_{K}^{l}(G)'$.  Trivially, both $\al^{a}, \bt_{a}$ are multiplicative.    
The map $J$ in (\ref{eq:btalt}) is of course just $\de^{-1}\circ \ga$.
\end{proof} 

We noted above that $\Ga$ is naturally a group.  We now turn to the natural, 
two-sided, jointly continuous action of the group $\Ga$ on $A(G)$. This is defined as follows:   for $x\in G$, define $xa,ax\in G$ by setting $xa=xa^{s(x)}, ax=a_{r(x)}x$. The continuity of this action follows from the continuity of the maps $u\to a^{u}, u\to a_{u}$ and that of the multiplication of $G$.
For $f\in C_{c}(G)$, define $af$, $fa\in C_{c}(G)$ by setting:
$af(x)=f(xa)$, $fa(x)=f(ax)$.  

\begin{proposition}  \label{prop:fgT}
Let $a\in \Ga, f,g\in C_{c}(G)$ and $\phi\in B(G)$, $T\in VN(G)$.  
Then $a(f*g)=f*ag$, $(f*g)a=fa*g$. Further, $a\phi,\, 
\phi a\in B(G)$, and these two functions have $B(G)$-norms equal to that of $\phi$.  Lastly, if $\phi\in \mathcal{A}(G)$, then $a\phi, \phi a\in \mathcal{A}(G)$, the norms of $\phi, a\phi, \phi a$ are all equal in $\mathcal{A}(G)$ and
$T(a\phi)=aT\phi$.
\end{proposition}
\begin{proof}
We have 
\begin{align*}
(f*g)a(x)& =\int f(t)g((a^{-1}t)^{-1}x)\,d\la^{r(ax)}(t)\\
&=\int f(az)g(z^{-1}x)\,d\la^{r(x)}(z)\\
=(fa*g)(x). 
\end{align*} 
Similarly, $a(f*g)=f*ag$.  
Next, write $\phi$ as a coefficient $(\xi,\eta)$.  Then $\phi a=(\xi,\eta')$ where
$\eta'(u)=(L(a_{u}))^{-1}\eta(J^{-1}u)$.  Then $\norm{\eta'}=\norm{\eta}$, and $\phi a\in B(G)$, $\norm{\phi a}=\norm{\phi}$.  
Similarly, $a\phi\in B(G)$ and $\norm{a\phi}=\norm{\phi}$. 
Next, let $\phi\in \mathcal{A}(G)$.  It is left to the reader to check that $a\phi, \phi a\in \mathcal{A}(G)$, and that $\norm{\phi}_{1}=\norm{a\phi}_{1}=\norm{\phi a}_{1}$. 
Lastly, by Proposition~\ref{prop:cont},
$T(a(f*g))=R_{ag}T(f)=T(f)*ag=a(T(f)*g)=aT(f*g)$.
\end{proof}

For our main theorem, we need to restrict the class of groupoids under consideration to those for which $r$ is locally trivial.

\vspace{.2in}
\noindent
{\bf Definition}\\  
A locally compact groupoid $G$ is said to be {\em locally a product} if the following holds:  for each $x_{0}\in G$, there exists an open neighborhood $U$ of $x_{0}$ in $G$, a locally compact Hausdorff space $Y$, a positive regular Borel measure $\mu$ on $Y$ and a homeomorphism $\Phi$ from $U$ onto $r(U)\x Y$ such that:
\bi
\item[(i)] $p_{1}(\Phi(x))=r(x)$ for all $x\in U$, where $p_{1}$ is the projection from 
$r(U)\x Y$ onto the first coordinate;
\item[(ii)] for each $u\in r(U)$ and with $\Phi^{u}$ the restriction of $\Phi$ to $U^{u}$, we have $(\Phi^{u})^{*}\mu=\la^{u}_{\mid U^{u}}$.
\ei
\vspace*{.2in}

Such an open set $U$ is called a {\em product open subset} of $G$. Examples of groupoids $G$ that are locally a product include Lie groupoids (more generally, continuous family groupoids (\cite{patcont})) and r-discrete groupoids.  

We now come to our main theorem.  We require two conditions on our groupoid $G$.  The first (i) of these is that $G$ is locally a product, i.e. is ``locally trivial'', has ``lots of'' local sections.  The second (ii) says that every point of $G$ lies on a global bisection.
(i) and (ii) are reasonable section conditions to require given that our duality theorem is formulated in terms of the group $\Ga$ of global {\em bisections}.  There are many examples of groupoids satisfying (i) and (ii).  For example, 
every ample groupoid (\cite[p.48]{patcont} does, and also the tangent groupoid for a manifold.

Under these conditions, the theorem says that the map $a\to (\al^{a},\bt_{a})$ is a homeomorphism from $\Ga$ onto $\Phi_{A(G)}$.  (In particular, $\Phi_{A(G)}$ is a group.)  It is easy to see that the map $a\to (\al^{a},\bt_{a})$ is one-to-one.  What requires more work (as it also does in the original group case considered by Eymard) is to show that the map is onto.

We now specify the topologies that we will use on $\Ga$ and $\Phi_{A(G)}$.  We will call each of these the {\em pointwise} topology.  Precisely, each $a\in \Ga$ is regarded as a function
$a:G^{0}\to G^{2}$, where $a(u)=(a_{u},a^{u})$.  We then give $\Ga$ the topology of pointwise convergence on $G^{0}$.  So $a_{n}\to a$ in $\Ga$ if and only if $a_{n}(u)\to a(u)$ in $G^{2}$ for all $u\in G^{0}$.  Turning to $\Phi_{A(G)}$, regard each $(\al,\bt)\in \Phi_{A(G)}$ as a function $(\al,\bt):A(G)\x G^{0}\to \C^{2}$ by: 
$(\al,\bt)(\phi,u)=(\al(\phi)(u),\bt(\phi)(u))$.  The pointwise topology on $\Phi_{A(G)}$ is then the topology of pointwise convergence on $A(G)\x G^{0}$.

\begin{theorem}        \label{th:main}
Assume 
\bi
\item[(i)] $G$ is locally a product.
\item[(ii)] if $x\in G$, then there exists $a\in \Ga$ such that $x\in a$.
\ei  
Then the map $\zeta$ taking $a\to (\al^{a},\bt_{a})$ is a homeomorphism for the pointwise topologies from $\Ga$ to $\Phi_{A(G)}$.
\end{theorem}
\begin{proof}  
The proof is an adaptation to the groupoid case of Eymard's proof that $G$ is the character space of $A(G)$ in the locally compact group case.
Let $(\al,\bt)\in \Phi_{A(G)}$ and $a\in \Ga$ be the element determined by $(\al,\bt)$.   
We want to show that $(\al,\bt)=(\al^{a},\bt_{a})$.  (This will give that $\zeta$ is onto.)
We show then that $\al=\al^{a}$, the result that $\bt=\al_{a}$ following using $G(r)$.
By Proposition~\ref{prop:mcag} and Proposition~\ref{prop:agba}, $\al$ determines a bounded linear map, also denoted by $\al$, from $\mathcal{A}(G)$ into $D$.
Define $\al a^{-1}:\mathcal{A}(G)\to D$ by:
$\al a^{-1}(\phi)=\al(a^{-1}\phi)$.  By Proposition~\ref{prop:fgT},
the map $\al a^{-1}$ is a bounded right module map.  Further since 
$f(\al a^{-1})=(a^{-1}f)\al$, it follows that 
$\al a^{-1}\in \mathcal{A}_{K}^{r}(G)'$.  Next, if $U$ is a neighborhood of $u\in G^{0}$ in $G$, then there exists $\phi\in A_{cf}(Ua)$ such that 
$\al(\phi)(a^{u})\ne 0$.  But then  $a\phi\in A_{cf}(U)$ 
and $\al a^{-1}(a\phi)(u)=\al(\phi)(u)\ne 0$.
So $G^{0}$ is the element of $\Ga$ determined by $\al a^{-1}$.  Clearly, 
$\al a^{-1}=\al^{G^{0}}$ if and only if 
$\al=\al^{a}$.  For the purposes of the theorem, we can therefore suppose that 
$a=G^{0}$.  

Let $T$ be the operator in $VN(G)$ determined by $\al_{\mid \mathcal{A}(G)}$.  So $\al_{T}=\al$ (Theorem~\ref{th:agdual}).  Then using (\ref{eq:alT}), $\al=\al^{G^{0}}$ if and only if
for all $f,g\in C_{c}(G)$, we have $\lan Tf,g\ran(u)=g*f^{*}(u)=\lan f,g\ran(u)$, i.e. if and only if 
$T$ is the identity map $I$.  So we have to show that $T=I$. 

We first show that  
$S(T\phi)\subset S(\phi)$ for $\phi\in A_{sp}(G)$.  For let $x\notin S(\phi)$ and $u=r(x)$.  By (ii), there exists $c\in \Ga$ such $x=c^{u}$.  Then by Proposition~\ref{prop:fgT}
and (\ref{eq:alT}),   
$T(\phi)(x)=(cT(\phi))(u)=\ov{\al((c\phi)^{*})(u)}$.
Now  $(c\phi)^{*}(u)=\ov{\phi(x)}$ and since $\phi=0$ on a neighborhood of $x$, 
it follows by continuity that $(c\phi)^{*}=0$ on a neighborhood of $u$ in $G$.  By Proposition~\ref{prop:AN}, (ii), (iv), it follows that 
$\al((c\phi)^{*})(u)=0$.  So $T(\phi)(x)=0$, and $T\phi$ vanishes on the complement of $S(\phi)$.  So $S(T\phi)\subset S(\phi)$.  

Now let $\Om$ be an open relatively compact subset of $G$.  Suppose that $\phi\in A_{sp}(G)$ is such that for each $u\in G^{0}$, the restriction $\phi_{\mid \Om^{u}}$ is constant, say $j_{u}$.  We will say that $\phi$ is {\em fiber constant} on $\Om$.  Note that the map $j$, where $j(u)=j_{u}$, is a continuous bounded function on $r(\Om)$. 

We will now show that $T\phi$ is also fiber constant on $\Om$.  Let $u\in r(\Om)$, and $q,p\in \Om^{u}$. We have to show that $T\phi(q)=T\phi(p)$.  To this end,
let $c,d\in \Ga$ be such that $c^{u}=q, d^{u}=p$, and $V$ be an open neighborhood of $u$ in $G$
such that $Vc\cup Vd\subset \Om$.  Let $U=Vc$ and $x\in U$.  Then $x,xc^{-1}d\in \Om^{r(x)}$.    Let  $\psi=(\phi - c^{-1}d\phi)\in A_{sp}(G)$.  Then $\psi$ is zero on $U$.  Since
$S(T\psi)\subset S(\psi)$, it follows that $T\psi$ is zero on $U$.  In particular, since $q\in U$, $0=T\psi(q)=(T\phi - T(c^{-1}d\phi))(q)=T\phi(q)-T\phi(p)$ as claimed.

Now let $f\in C_{c}(G)$ and suppose 
that the support $S(f)$ of $f$ lies in a relatively compact, product open subset $U$.  In the preceding notation, we will identify $U$ with a product $r(U)\x Y$, with associated r-fiber preserving homeomorphism $\Phi:U\to r(U)\x Y$ and measure $\mu$ on $Y$.  Let 
$L$ be a compact subset of $Y$ for which $\mu(L\setminus L^{o})=0$.  Let $V_{n}$ be a sequence of open subsets of $Y$ with $V_{n}\subset \ov{V_{n}}\subset V_{n+1}$, $\cup V_{n}=L^{o}$ and such that $\mu(L\setminus \ov{V_{n}})=\mu(L^{o}\setminus \ov{V_{n}})<1/n$.  Let $Z$ be an open subset of $r(U)$ such that $\ov{Z}$ is compact and contained in $r(U)$.    
Let $b\in C_{0}(Z)$.
By Proposition~\ref{prop:eleag2}, (ii),
there exists $\phi_{n}\in A_{cf}(r(U)\x L^{o})$ such that $0\leq \phi_{n}\leq 1$ and $\phi_{n}(z,y)=1$ for all $(z,y)$ in a neighborhood $V$ of $\ov{Z}\x \ov{V_{n}}$.  Let 
$F(z,x)=b(z)\chi_{L}(x)$ for $(z,x)\in r(U)\x Y$.  For $z\in Z$, we have  
$\norm{\phi_{n}b-F}^{z}\leq \norm{b}_{\infty}(1/n)^{1/2}$.
So $\norm{\phi_{n}b-F}\to 0$, in $E^{2}$, and by continuity, 
\begin{equation}       \label{eq:tbphi}
T(\phi_{n}b)\to TF \text{ in } E^{2}.
\end{equation}
Next, $\phi_{n}b$ is fiber constant on $V$, and so $T(\phi_{n}b)$
also has constant fiber on $V$.  So for some continuous function $k_{n}$ on $r(V)$, we have $T(\phi_{n}b)(z,x)=k_{n}(z)$ on $\ov{Z}\x V_{n}$.  Further $k_{n}$ is independent of $n$
since $(\phi_{n+1}b)_{\mid Z\x V_{n}}=(\phi_{n}b)_{\mid Z\x V_{n}}$.  Next, since 
$S(T(\phi_{n}b))\subset S(\phi_{n}b)$,  we have $k_{n}$ vanishing on $r(V)\setminus \ov{Z}$, so that $k_{n}\in C_{0}(Z)$.  Also $S(k_{n})\subset S(b)$.  Write
$k_{b}$ in place of $k_{n}$.  
Note that $TF=\chi_{Z\x L}k_{b}$ using (\ref{eq:tbphi}).  

Let $F'\in C_{c}(r(U)\x Y)$ be such that $0\leq F'\leq 1$ and $F'=1$ on $\ov{Z}\x L$.
Then for $z\in Z$ and any $n$, we have 
\begin{align*}
\mid k_{b}(z)\mid\mu(L)^{1/2} & =\lim_{n}\norm{T(\phi_{n}b)}^{z} \\
&\leq \norm{T}\norm{b}_{\infty}\lim_{n}\norm{\phi_{n}}^{z} \\
&\leq \norm{T}\norm{b}_{\infty}\norm{F'}. 
\end{align*}
It follows that the map $R:C_{0}(Z)\to C_{0}(Z)$, where $Rb=k_{b}$,
is a linear, bounded map.   Further, $S(Rb)\subset S(b)$ for all $b$.  
So $R$ satisfies the conditions of Proposition~\ref{prop:supp}, and there exists a bounded continuous function $k$ on $Z$ such that $T(f)=fk$ for $f$ of the form $\chi_{Z\x L}b$.

Then  contract down onto a general compact subset $K$ of $Y$ by open relatively compact sets $L$ with null boundary.  We get the same $k$ for each $L$,  and obtain that each function
$\chi_{Z\x K}b\in E^{2}$ and $T(\chi_{Z\x K}b)=(k\circ r)(\chi_{Z\x K}b)$.  Next
for $g\in C_{c}(Y)$, approximate in $L^{2}(\mu)$ the function $g$ by linear combinations of 
$\chi_{K}$'s to get $T(b\otimes g)=(k\circ r)(b\otimes g)$.  Using the Stone-Weierstrass theorem, we then get $T(f)=fk$ for all $f\in C_{c}(Z\x Y)$.   

Let $U'$ be the product open set in $G$ corresponding to $Z\x Y$.  
The functions $k$ for different choices of such $U'$ are all compatible. 
So there exists a continuous function $h$ such that whenever $U'$ is a such a product open subset of $G$ and $f\in C_{c}(U')$, then $Tf=fh$.  The latter equality is true for all $f\in C_{c}(G)$ by a partition of unity argument.  So for all $f\in C_{c}(G)$,
\begin{equation}  \label{eq:thrf}
Tf=fh.
\end{equation}
It follows from (\ref{eq:thrf}) that $\norm{h}_{\infty}\leq \norm{T}$ so that $h\in C(G)$.

Now for $f,g\in C_{c}(G)$, we have 
\[ \al(g*f^{*})(u)=\lan Tf,g\ran(u)=\ov{h(u)}\lan f,g\ran(u) = \ov{h(u)}g*f^{*}(u).   \]
So $\al(\phi)(u)=\ov{h(u)}\phi(u)$ for all $\phi\in A_{cf}(G)$, and since $\al$ is multiplicative, we get $\ov{h(u)}^{2}\phi(u)\phi_{1}(u)=\ov{h(u)}\phi(u)\phi_{1}(u)$ for all $\phi,\phi_{1}\in A_{cf}(G)$. It follows that $h(u)$ is either $1$ or $0$.  In fact $h(u)=1$ since $\al(A_{c}(G))(u)\ne \{0\}$ for all $u\in G^{0}$.  So $T=I$ and 
$\al=\al_{G^{0}}$ as we had to prove.

Using (\ref{eq:albta}) and the fact that $A(G)$ separates the points of $G$ 
(Proposition\ref{prop:eleag2}, (ii)) it follows that $\zeta$ is a homeomorphism.
\end{proof}

\vspace{.2in}
\noindent
\begin{center}
\bf{Some open problems}
\end{center}
All of the following questions are answered positively for $G_{n}$ and for all locally compact groups.
\bi
\item[(i)]  Is it true that $A_{cf}(G)=A_{sp}(G)=\mathcal{A}(G)=A(G)$?  
\item[(ii)] Is the character space of $A(G)$ equal to $G$?   
\item[(iii)] Is $\Phi^{r}_{A(G)}=\De^{r}$? 
\item[(iv)] When is $\norm{.}_{cb}=\norm{.}$ on $B(G)$?  
\item[(v)] If $G$ is amenable, does $A(G)$ have a bounded approximate identity, and is $B(G)$ the multiplier algebra of $A(G)$.
\ei


\end{document}